\documentclass[12pt]{amsart}
\usepackage{amsmath,amssymb}
\newtheorem{theorem}{Theorem}
\newtheorem{corollary}{Corollary}
\newtheorem{lemma}{Lemma}
\newtheorem{proposition}{Proposition}
\begin{document}
\centerline{\Large On the pseudohermitian sectional curvature}
\vskip 0.25cm \centerline{\Large of a strictly pseudoconvex CR
manifold} \vskip 0.5cm \centerline{Elisabetta {\large
Barletta}\footnote{Universit\`a degli Studi della Basilicata,
Dipartimento di Matematica, Contrada Macchia Romana, 85100
Potenza, Italy, e-mail: {\tt barletta@unibas.it} 
}}
\title{}
\author{}
\maketitle
\begin{abstract}
We show that the pseudohermitian sectional curvature $H_\theta
(\sigma )$ of a contact form $\theta$ on a strictly pseudoconvex
CR manifold $M$ measures the difference between the lengthes of a
circle in a plane tangent at a point of $M$ and its projection on
$M$ by the exponential map associated to the Tanaka-Webster
connection of $(M , \theta )$. Any Sasakian manifold $(M , \theta
)$ whose pseudohermitian sectional curvature $K_\theta (\sigma )$
is a point function is shown to be Tanaka-Webster flat, and hence
a Sasakian space form of $\varphi$-sectional curvature $c = -3$.
We show that the Lie algebra ${\mathfrak i}(M , \theta )$ of all
infinitesimal pseudohermitian transformations on a strictly
pseudoconvex CR manifold $M$ of CR dimension $n$ has dimension
$\leq (n+1)^2$ and if $\dim_{\mathbb R} {\mathfrak i}(M , \theta )
= (n+1)^2$ then $H_\theta (\sigma ) =$ constant.
\end{abstract}
\section{Introduction}
In his famous 1987 paper S. M. Webster introduced (cf.
\cite{kn:Web1}) the notion of pseudohermitian sectional curvature
$H_\theta$ of a nondegenerate CR manifold, associated to a fixed
contact form $\theta$, and exhibited a class of spherical
nondegenerate real hypersurfaces $M \subset {\mathbb C}^{n+1}$
with $H_\theta (\sigma ) = \pm 1/(2c)$, for each $c \in (0 , +
\infty )$. If $M$ is a nondegenerate CR manifold and $\theta$ a
contact form on $M$ then let $R$ be the curvature $4$-tensor field
of the Tanaka-Webster connection $\nabla$ of $(M, \theta )$. Let
$\sigma \subset H(M)_x$ be a {\em holomorphic} $2$-{\em plane}
tangent at $x \in M$ i.e. $J_x (\sigma ) = \sigma$. Here $H(M)$ is
the maximal complex distribution of $M$ and $J$ its complex
structure. If $\{ X , J_x X \}$ is a linear basis of $\sigma$ then
we set
\begin{equation}
\label{e:1} H_\theta (\sigma ) = \frac{1}{4} G_{\theta ,
x}(X,X)^{-2} R_x (X , J_x X , X , J_x X ) \;.
\end{equation}
The definition of $H_\theta (\sigma )$ doesn't depend upon the
choice of basis in $\sigma$ because of $R(X,Y,Z,W) = - R(X,Y,W,Z)$
(as the curvature is a $2$-form) and $R(X,Y,Z,W) = - R(Y,X, Z,W)$
(as the Levi form is parallel with respect to $\nabla$). Then
$H_\theta$ is a ${\mathbb R}$-valued function on the total space
of the Grassmann bundle $G_2 ({\mathbb C}^n ) \to G_2 (H(M))
\stackrel{\pi}{\longrightarrow} M$ of all holomorphic $2$-planes
tangent to $M$. We also set $H_\theta = {\rm Sect}(M, \theta )$.
The coefficient $1/4$ in \eqref{e:1} is chosen such that the
standard sphere $S^{2n+1} \subset {\mathbb C}^{n+1}$ together with
the canonical contact form $\theta_0 = \frac{i}{2}
(\overline{\partial} -
\partial )|z|^2$ has constant pseudohermitian sectional curvature
${\rm Sect}(S^{2n+1} , \theta_0 )
 \equiv 1$. Clearly \eqref{e:1} is a (pseudohermitian) analog of the
{\em holomorphic} sectional curvature of a Hermitian manifold (cf.
e.g. \cite{kn:KoNo}, Vol. II, p. 168) rather than an analog of the
sectional curvature of a Riemannian manifold (cf. \cite{kn:KoNo},
Vol. I, p. 202). Yet if $G_2 ({\mathbb R}^{2n+1}) \to G_2 (T(M))
\stackrel{\pi}{\longrightarrow} M$ is the Grassman bundle of all
$2$-planes tangent at $M$ then \eqref{e:1} is the restriction to
$G_2 (H(M))$ of the function $K_\theta : G_2 (T(M)) \to {\mathbb
R}$ given by
\begin{equation}
\label{e:2} K_\theta (\sigma ) = R_x (u , v, u , v), \;\;\; \sigma
\subset T_x (M) \; ,
\end{equation}
where $\{ u , v \}$ is a $g_{\theta , x}$-orthonormal basis of
$\sigma$ and $g_\theta$ is the Webster metric of $(M , \theta )$
(cf. Section 2 for definitions) and \eqref{e:2} may be referred to
as the ({\em pseudohermitian}) {\em sectional curvature}
determined by the (arbitrary) $2$-plane $\sigma$.
\par
A number of fundamental questions remain unanswered. First, what
is the geometric interpretation of $K_\theta (\sigma )$?
Precisely, if $\sigma \in G_2 (T(M))_x$ and $r \, w(s) = r (\cos
s) u + r (\sin s) v$ is a circle in $\sigma$ and $\beta_r(s) =
\exp_x r w(s)$, $0 \leq s \leq 2\pi$, then is $K_\theta (\sigma )$
a ``measure'' of the difference $2\pi r - L(\beta_r)$? Here
$\exp_x$ is the exponential map associated to the Tanaka-Webster
connection $\nabla$ of $(M , \theta )$ and $L(\beta_r)$ the length
of $\beta_r$. Another fundamental question is whether the
algebraic machinery in \cite{kn:KoNo} (cf. Vol. I, p. 198-203, and
Vol. II, p. 165-169) applies, eventually leading to a meaningful
concept of pseudohermitian space form. Moreover, as
pseudohermitian transformations are (within pseudohermitian
geometry) analogs to isometries between Riemannian manifolds, it
is a natural question whether manifolds $(M , \theta )$ whose Lie
algebra ${\mathfrak i}(M , \theta )$ of infinitesimal
pseudohermitian transformations has maximal dimension have
constant pseudohermitian sectional curvature.
\par
Our findings are that the pseudohermitian sectional curvature
\eqref{e:1} satisfies
\begin{equation}
L(\beta_r) = 2\pi r - \frac{\pi r^3}{12} (16 H_\theta (\sigma ) -
3) + O(r^4 ) \label{e:3}
\end{equation}
(cf. Theorem \ref{t:1} below for the precise statement) providing
the geometric interpretation mentioned above. Also we prove a
Schur like result, cf. Theorem \ref{t:2} below. Combining Theorem
\ref{t:2} with a result by Y. Kamishima, \cite{kn:Kam}, we obtain
\begin{corollary} Let $(M, \theta )$ be a compact connected Sasakian manifold of CR
dimension $n \geq 2$. If there is a $C^\infty$ function $f : M \to
{\mathbb R}$ such that $K_\theta = f \circ \pi$ then $M$ is
isometric to the Heinsenberg infra\-nil\-ma\-ni\-fold ${\mathbb
H}_n /\Gamma$ {\rm (}with $\Gamma = \rho (\pi_1 (M)) \subset
{\mathbb H}_n \rtimes {\rm U}(n)${\rm )}. \label{c:1}
\end{corollary}
Here ${\mathbb H}_n$ is the Heisenberg group endowed with the
standard strictly pseudoconvex CR structure and canonical contact
form (cf. e.g. \cite{kn:DrTo}, Chapter 1).
\par
The paper is organized as follows. Section 1 is devoted to a
remainder of CR and pseudohermitian geometry and to the proof of
Theorem \ref{t:1}. The main technical ingredient are Jacobi fields
of the Tanaka-Webster connection, on the line of thought in
\cite{kn:BaDr}. A Schur like result for the sectional curvature
\eqref{e:2} and the proof of Corollary \ref{c:1} form the object
of Section 3. In Section 4 we show (cf. Theorem \ref{t:3} below)
that for any strictly pseudoconvex CR manifold $\dim_{\mathbb R}
{\mathfrak i}(M , \theta ) \leq (n+1)^2$ and if $\dim_{\mathbb R}
{\mathfrak i}(M , \theta ) = (n+1)^2$ then $(M , \theta )$ has a
constant pseudohermitian sectional curvature \eqref{e:1}. The
proof of Theorem \ref{t:3} relies on standard techniques in the
theory of (infinitesimal) affine transformations. The explicit
expression of the curvature tensor of a pseudohermitian space form
(i.e. a pseudohermitian manifold whose sectional curvature
\eqref{e:1} is constant) is derived in Section 5 (cf. \eqref{e:19}
in Theorem \ref{t:4} below) paving the road towards a study of the
geometry of the second fundamental form of a CR submanifold of a
pseudohermitian space form (in the spirit of \cite{kn:YaKo}, p.
76-136). The computational details (leading, as a byproduct, to a
Sasakian version of the K\"ahlerian Schur theorem) are provided in
the Appendix A to this paper. A classification result of E. Musso,
\cite{kn:Mus}, and our Theorem \ref{t:4} lead to
\begin{corollary} Let $(M, \theta )$ be a $G$-homogeneous
pseudohermitian space form of pseudohermitian sectional curvature
$H_\theta (\sigma ) = c$, $c \in {\mathbb R}$, with $L_\theta$
positive definite. {\rm i)} If $c > 0$ then $(M, \theta )$ is
contact homothetic to the canonical pseudohermitian manifold of
index $k$ over $B$. {\rm ii)} If $c < 0$ then $(M, \theta )$ is
contact homothetic to either $B \times S^1$ or $B \times {\mathbb
R}$. {\rm iii)} If $c = 0$ then $(M , \theta )$ is contact
homothetic to either ${\mathbb C}^n \times S^1$ or ${\mathbb C}^n
\times {\mathbb R}$. \label{c:2}
\end{corollary}
The description of the pseudohermitian structures on the model
spaces i)-iii) in Corollary \ref{c:2} is provided in Section 5.
Finally, in Section 6 we show that given a pseudohermitian
immersion $f : M \to M^\prime$ between two strictly pseudoconvex
CR manifolds the sectional curvature \eqref{e:1} of $M$ doesn't
exceed the sectional curvature \eqref{e:1} of the ambient space.
Theorem \ref{t:6} in Section 6 is suitable for several
applications. For instance
\begin{corollary}
There is no pseudohermitian immersion of the standard sphere
$S^{2m+1}$ into an ellipsoid $\{ (z, w) :
g_{\alpha\overline{\beta}} z^\alpha \overline{z}^\beta - w
\overline{w} + c = 0\} \subset {\mathbb C}^{n+1}$, with $c \in (0,
+ \infty )$ and $[g_{\alpha \overline{\beta}}] \in \mathrm{GL}(n,
\mathbb{C})$ Hermitian. \label{c:3}
\end{corollary}
\begin{corollary} For any compact Sasakian manifold $(M, \theta )$ there are $n \geq 1$ and $A = \{ 0 < a_1 \leq a_2 \leq \cdots
\leq a_{n+1} \}$ such that $H_\theta (\sigma ) \leq {\rm
Sect}(S^{2n+1} , \theta_A )$ where $\theta_A = \left(
\sum_{j=1}^{n+1} a_j |z_j |^2 \right)^{-1} \theta_0$. \label{c:4}
\end{corollary}
\section{The geometric interpretation of pseudohermitian sectional
curvature}
\subsection{The Tanaka-Webster connection} Let us start  by recalling the notions of CR and
pseudohermitian geometry needed through this paper. Let $(M,
T_{1,0}(M))$ be a $(2n+d)$-dimensional CR manifold, of CR
dimension $n$, where $T_{1,0}(M)$ denotes the CR structure. The
maximal complex distribution is $H(M) = {\rm Re}\{ T_{1,0}(M)
\oplus T_{0,1}(M) \}$. It carries the complex structure $J : H(M)
\to H(M)$ given by $J(Z + \overline{Z}) = i(Z - \overline{Z})$,
for any $Z \in T_{1,0}(M)$. Throughout $T_{0,1}(M) =
\overline{T_{1,0}(M)}$ and overlines denote complex conjugates.
The standard example of a CR manifold is that of a real
submanifold $M \subset {\mathbb C}^N$ such that $\dim_{\mathbb C}
[T_x (M) \otimes_{\mathbb R} {\mathbb C}] \cap T^{1,0}({\mathbb
C}^N )_x =$ constant, $x \in M$. This is of course always true for
real hypersurfces in ${\mathbb C}^N$.
\par
On each CR manifold $M$ there is a natural first order
differential operator $\overline{\partial}_b$ given by
$(\overline{\partial}_b f) \overline{Z} = \overline{Z}(f)$ for any
$C^1$ function $f : M \to {\mathbb C}$ and any $Z \in T_{1,0}(M)$.
Then $\overline{\partial}_b f = 0$ are the {\em tangential
Cauchy-Riemann equations} and a $C^1$ solution is a {\em CR
function} on $M$.
\par
Let $H(M)^\bot \subset T^* (M)$ be the conormal bundle associated
to $H(M)$. When $M$ has hypersurface type (i.e. $d = 1$) and $M$
is orientable, which we shall always assume, $H(M)^\bot$ is a
trivial line bundle hence $M$ admits globally defined nowhere zero
differential $1$-forms $\theta$ such that ${\rm Ker}(\theta ) =
H(M)$. These are referred to as {\em pseudohermitian structures}.
With each pseudohermitian structure $\theta$ one may associate the
Levi form $L_\theta (Z, \overline{W}) = - i (d \theta )(Z,
\overline{W})$, $Z, W \in T_{1,0}(M)$, and $M$ is {\em
nondegenerate} (respectively {\em strictly pseudoconvex}) if
$L_\theta$ is nondegenerate (respectively positive definite) for
some $\theta$. Two pseudohermitian structures $\theta$ and
$\hat{\theta}$ are related by $\hat{\theta} = f \; \theta$ for
some $C^\infty$ function $f : M \to {\mathbb R} \setminus \{ 0 \}$
and a simple calculation shows that $L_{\hat{\theta}} = f
L_\theta$. Nondegeneracy is a {\em CR invariant} property i.e. it
is invariant under a transformation $\hat{\theta} = f \theta$.
Clearly, strict pseudoconvexity is not a CR invariant property
(e.g. if $L_\theta$ is positive definite and $\hat{\theta} = -
\theta$ then $L_{\hat{\theta}}$ is negative definite). If $M$ is a
nondegenerate CR manifold of CR dimension $n$ then each
pseudohermitian structure is a {\em contact form} i.e. $\theta
\wedge (d \theta )^n$ is a volume form on $M$.
\par
Let $M$ be a nondegenerate CR manifold and $\theta$ a contact form
on $M$. The pair $(M , \theta )$ is commonly referred to as a {\em
pseudohermitian manifold}. There is a unique nowhere zero globally
defined tangent vector field $T$ on $M$, transverse to $H(M)$,
determined by $\theta (T) = 1$ and $(d \theta )(T, X) = 0$ for any
$X \in T(M)$ ($T$ is the {\em characteristic direction} of $d
\theta$). On any pseudohermitian manifold $(M, \theta )$ there is
a unique linear connection $\nabla$ (the {\em Tanaka-Webster
connection} of $(M , \theta )$) such that i) $H(M)$ is parallel
with respect to $\nabla$, ii) $\nabla J = 0$ and $\nabla g_\theta
= 0$, and iii) the torsion $T_\nabla$ of $\nabla$ is {\em pure}
i.e.
\[ T_\nabla (Z, W) = 0, \;\;\; T_\nabla (Z , \overline{W}) = 2 i
L_\theta (Z, \overline{W}) T, \;\;\; Z,W \in T_{1,0}(M) \; , \]
\[ \tau \circ J + J \circ \tau = 0 \; . \]
Cf. N. Tanaka, \cite{kn:Tan}, S. M. Webster, \cite{kn:Web1}, or
Chapter I of \cite{kn:DrTo}. Here $g_\theta$ is the {\em Webster
metric} i.e. the semi-Riemannian metric on $M$ defined by
\[ g_\theta (X,Y) = (d \theta )(X, J Y), \;\;\; g_\theta (T , X) =
0, \;\;\; g_\theta (T,T) = 1 \; , \] for any $X,Y \in H(M)$. Also
$\tau$ is the {\em pseudohermitian torsion} i.e. the vector-valued
$1$-form $\tau (X) = T_\nabla (T, X)$, $X \in T(M)$. The complex
structure $J : H(M) \to H(M)$ appearing in axiom (ii) is thought
of as extended to a $(1,1)$-tensor field on $M$ by requesting that
$J T = 0$. When $M$ is strictly pseudoconvex and $L_\theta$ is
positive definite the Webster metric is a Riemannian metric on $M$
and $(J, T, \theta , g_\theta )$ is a contact metric structure (in
the sense of D. E. Blair, \cite{kn:Bla}) which is normal if and
only if $\tau = 0$. If this is the case then $g_\theta$ is a
Sasakian metric on $M$. Therefore Sasakian manifolds are precisely
the strictly pseudoconvex CR manifolds with a fixed contact form
$\theta$ such that the Levi form $L_\theta$ is positive definite
and the pseudohermitian torsion of the Tanaka-Webster connection
vanishes. By a result of G. Marinescu et al., \cite{kn:MaYe}, for
any Sasakian manifold $M$ there is a CR embedding $M \to {\mathbb
C}^N$ for some $N \geq 2$.

\subsection{Jacobi fields} A study of Jacobi fields of the
Tanaka-Webster connection on a nondegenerate CR manifold was
started in \cite{kn:BaDr}. Let $M$ be a strictly pseudoconvex CR
manifold, of CR dimension $n$, and $\theta$ a contact form with
$L_\theta$ positive definite. Let $x \in M$ and let $\exp_x$ be
the exponential mapping, associated to the Tanaka-Webster
connection $\nabla$ of $(M , \theta )$. Here we use a few facts
from the general theory of linear connections on manifolds e.g. by
Proposition 8.2 in \cite{kn:KoNo}, Vol. I, p. 147, there is $r_0 >
0$ such that $\exp_x : B(x, r_0 ) \to M$ is a $C^\infty$
diffeomorphism on some neighborhood $U$ of $x$ in $M$. Here $B(x ,
r_0 ) = \{ v \in T_x (M) : \| v \| < r_0 \}$ and $\| v \|^2 =
g_{\theta , x} (v,v)$. Let
\[ L(\beta_r) = \int_0^{2\pi} \| \dot{\beta}_r (s) \| \; d s \]
be the length of the curve $\beta_r$ (defined in the Introduction)
in $(M , g_\theta )$. Let $\gamma_v (t) = \exp_x t v$ denote the
geodesic of $\nabla$ of initial conditions ($x, v)$, $v \in T_x
(M)$. Given $0 < r < r_0$ we consider the geodesics $\gamma_{w(s)}
: [-r , r] \to U$ and set $\beta_t (s) = \gamma_{w(s)}(t)$. Next
let $X_s$ be the vector field along $\gamma_{w(s)}$ defined by
\[ X_{s, \gamma_{w(s)}(t)} = \dot{\beta}_t (s), \;\;\; 0 \leq s
\leq 2\pi , \;\;\; |t| \leq r \; . \] Then $L(\beta_r ) =
\int_0^{2\pi} \| X_s \|_{\gamma_{w(s)} (r)} \, d s$. Once again a
general fact within connection theory (cf. Theorem 1.2 in
\cite{kn:KoNo}, Vol. II, p. 64) guarantees that $X_s$ is a Jacobi
field of the Tanaka-Webster connection i.e. $X_s$ satisfies the
Jacobi equation
\begin{equation} \nabla^2_{\dot{\gamma}_{w(s)}} X_s +
\nabla_{\dot{\gamma}_{w(s)}} T_\nabla (X_s , \dot{\gamma}_{w(s)} )
+ R(X_s , \dot{\gamma}_{w(s)}) \dot{\gamma}_{w(s)} = 0 \label{e:4}
\end{equation}
along $\gamma_{w(s)}$. Let us set $X^\prime_s =
\nabla_{\dot{\gamma}_{w(s)}} X_s$ for simplicity. An elementary
calculation shows that $X_s$ satisfies the initial conditions
\begin{equation}
X_{s,x} = 0, \;\;\; X^\prime_{s,x} = w(s + \frac{\pi}{2}) \;.
\label{e:5}
\end{equation}
We wish to write the Taylor development of $f(r) = \| X_s
\|^2_{\gamma_{w(s)}(r)}$ (with $0 \leq s \leq 2 \pi$ fixed) up to
order $4$. This is the classical approach to the geometric
interpretation of sectional curvature in Riemannian geometry,
except that we must deal with the presence of torsion terms. The
first of the initial conditions \eqref{e:5} gives $f(0) = 0$.
Next, as $\nabla g_\theta = 0$
\begin{equation}
f^\prime (r) = 2 g_\theta (X_s^\prime , X_s )_{\gamma_{w(s)}(r)}
\label{e:6}
\end{equation}
hence $f^\prime (0) = 0$. Differentiating in \eqref{e:6} we obtain
\begin{equation}
f^{\prime\prime}(r) = 2 g_\theta (\nabla_{\dot{\gamma}_{w(s)}}^2
X_s , X_s )_{\gamma_{w(s)}(r)} + 2 \| X_s^\prime
\|^2_{\gamma_{w(s)}(r)} \label{e:7}
\end{equation}
hence (by \eqref{e:5})
\[ f^{\prime\prime}(0) = 2 \| w(s + \frac{\pi}{2}) \|^2 = 2 \; . \]
Let us set $P_s = \nabla_{\dot{\gamma}_{w(s)}} T_\nabla$ for
simplicity. Similarly we may differentiate in \eqref{e:7} so that
to get
\begin{equation}
f^{\prime\prime\prime}(r) = 2 g_\theta
(\nabla_{\dot{\gamma}_{w(s)}}^3 X_s , X_s )_{\gamma_{w(s)}(r)} + 6
g_\theta (\nabla^2_{\dot{\gamma}_{w(s)}} X_s , X^\prime_s
)_{\gamma_{w(s)}(r)} \label{e:8}
\end{equation}
hence (by the Jacobi equation \eqref{e:4})
\[
f^{\prime\prime\prime}(0) = 6 g_\theta
(\nabla^2_{\dot{\gamma}_{w(s)}} X_s , X^\prime_s )_x = \]
\[ =
- 6 g_\theta (\nabla_{\dot{\gamma}_{w(s)}} T_\nabla (X_s ,
\dot{\gamma}_{w(s)}) , X^\prime_s )_x - 6 g_\theta (R(X_s ,
\dot{\gamma}_{w(s)} ) \dot{\gamma}_{w(s)} , X^\prime_s )_x =
\]
\[ = - 6 g_\theta (P_s (X_s , \dot{\gamma}_{w(s)}),
X^\prime_s )_x + g_\theta ( T_\nabla (X^\prime_s ,
\dot{\gamma}_{w(s)}), X_s^\prime )_x = \]
\[= -
6 \langle T_{\nabla , x} (w(s+ \frac{\pi}{2}), w(s)) \, , \, w(s +
\frac{\pi}{2}) \rangle \] where $g_{\theta , x} = \langle \; , \;
\rangle$. Thus
\[ f^{\prime\prime\prime}(0) = 6 \langle T_{\nabla , x}(u,v) \, ,
\, w(s + \frac{\pi}{2}) \rangle \; . \] Finally we may
differentiate in \eqref{e:8} to obtain
\[
f^{(4)}(r) = 2 g_\theta (\nabla^4_{\dot{\gamma}_{w(s)}} X_s \, ,\,
X_s )_{\gamma_{w(s)}(r)} + \]
\[ + 8 g_\theta
(\nabla_{\dot{\gamma}_{w(s)}}^3 X_s \, , \, X_s^\prime
)_{\gamma_{w(s)}(r)} + 6 \| \nabla_{\dot{\gamma}_{w(s)}}^2 X_s
\|_{\gamma_{w(s)}(r)}^2  \; .\] Let us evaluate the terms in the
right hand side at $r = 0$. The first term vanishes (by
\eqref{e:5}). To compute the second term note first that (by
\eqref{e:4})
\[
\nabla^2_{\dot{\gamma}_{w(s)}} T_\nabla (X_s ,
\dot{\gamma}_{w(s)}) = \nabla_{\dot{\gamma}_{w(s)}} \{ P_s (X_s ,
\dot{\gamma}_{w(s)} ) + T_\nabla (X_s^\prime ,
\dot{\gamma}_{w(s)}) \} = \]
\[= (\nabla_{\dot{\gamma}_{w(s)}}
P_s )(X_s , \dot{\gamma}_{w(s)} ) + 2 P_s (X^\prime_s ,
\dot{\gamma}_{w(s)}) + T_\nabla (\nabla^2_{\dot{\gamma}_{w(s)}}
X_s , \dot{\gamma}_{w(s)}) = \]
\[ =
(\nabla_{\dot{\gamma}_{w(s)}} P_s )(X_s , \dot{\gamma}_{w(s)} ) +
2 P_s (X^\prime_s , \dot{\gamma}_{w(s)}) -
\]
\[
- T_\nabla (\nabla_{\dot{\gamma}_{w(s)}} T_\nabla (X_s ,
\dot{\gamma}_{w(s)}) , \dot{\gamma}_{w(s)}) - T_\nabla (R(X_s ,
\dot{\gamma}_{w(s)}) \dot{\gamma}_{w(s)} , \dot{\gamma}_{w(s)}) =
\]
\[
= (\nabla_{\dot{\gamma}_{w(s)}} P_s )(X_s , \dot{\gamma}_{w(s)} )
+ 2 P_s (X^\prime_s , \dot{\gamma}_{w(s)}) - T_\nabla (P_s (X_s ,
\dot{\gamma}_{w(s)}), \dot{\gamma}_{w(s)}) -
\]
\[
- T_\nabla (T_\nabla (X_s^\prime , \dot{\gamma}_{w(s)}) ,
\dot{\gamma}_{w(s)}) - T_\nabla (R(X_s ,
\dot{\gamma}_{w(s)})\dot{\gamma}_{w(s)} , \dot{\gamma}_{w(s)})
\] hence \begin{equation} (\nabla^2_{\dot{\gamma}_{w(s)}}
T_\nabla (X_s , \dot{\gamma}_{w(s)}))_x = T_{\nabla , x}
(T_{\nabla , x}(u,v) , w(s)) - 2 P_{s,x}(u,v) \; . \label{e:9}
\end{equation}
Similarly
\[ \nabla_{\dot{\gamma}_{w(s)}} R(X_s , \dot{\gamma}_{w(s)})
\dot{\gamma}_{w(s)} = (\nabla_{\dot{\gamma}_{w(s)}} R)(X_s ,
\dot{\gamma}_{w(s)}) \dot{\gamma}_{w(s)} + R(X^\prime_s ,
\dot{\gamma}_{w(s)}) \dot{\gamma}_{w(s)} \] hence
\begin{equation}
(\nabla_{\dot{\gamma}_{w(s)}} R(X_s , \dot{\gamma}_{w(s)})
\dot{\gamma}_{w(s)})_x = R_x (w(s+ \frac{\pi}{2}), w(s)) w(s) \; .
\label{e:10}
\end{equation}
Therefore (by \eqref{e:4} and \eqref{e:9}-\eqref{e:10})
\[
g_\theta (\nabla^3_{\dot{\gamma}_{w(s)}} X_s , X_s^\prime )_x =
\]
\[
= - g_{\theta} (\nabla_{\dot{\gamma}_{w(s)}} T_\nabla (X_s ,
\dot{\gamma}_{w(s)}), X_s^\prime )_x - g_\theta
(\nabla_{\dot{\gamma}_{w(s)}} R(X_s ,
\dot{\gamma}_{w(s)})\dot{\gamma}_{w(s)} , X_s^\prime )_x =
\]
\[ = 2 \langle P_{s,x}(u,v) , w(s +
\frac{\pi}{2}) \rangle - \langle T_{\nabla , x} (T_{\nabla ,
x}(u,v), w(s)), w(s + \frac{\pi}{2}) \rangle - \]
\[ - \langle R_x (w(s+\frac{\pi}{2}), w(s))w(s) \, , \,
w(s+\frac{\pi}{2})\rangle \; .\] Finally
\[
\| \nabla^2_{\dot{\gamma}_{w(s)}} X_s \|^2_x = \|
\nabla_{\dot{\gamma}_{w(s)}} T_\nabla (X_s , \dot{\gamma}_{w(s)})
+ R(X_s , \dot{\gamma}_{w(s)})\dot{\gamma}_{w(s)} \|^2_x =
\]
\[= \| P_s (X_s , \dot{\gamma}_{w(s)}) + T_\nabla
(X_s^\prime , \dot{\gamma}_{w(s)})\|_x^2 = \| T_{\nabla ,
x}(u,v)\|^2 \] and we may conclude that
\[
f^{(4)}(0) = 6 \| T_{\nabla , x}(u,v)\|^2 - 32 K_\theta (\sigma )
+ \]
\[ + 16 \langle P_{s,x}(u,v) , w(s+
\frac{\pi}{2}) \rangle - 8 \langle T_{\nabla , x} (T_{\nabla , x}
(u,v) , w(s)) \, , \, w(s + \frac{\pi}{2})\rangle \; .
\] We obtain the following
\begin{theorem} Let $M$ be a strictly pseudoconvex CR manifold and
$\theta$ a contact form on $M$ such that $L_\theta$ is positive
definite. Then
\begin{equation}
H_\theta (\sigma ) = \frac{3}{16} + \lim_{r \to 0} \frac{3}{4\pi
r^3} \{ 2 \pi r - L(\beta_r )\} \; . \label{e:11}
\end{equation}
for any holomorphic $2$-plane $\sigma \subset H(M)_x$ and $x \in
M$. \label{t:1}
\end{theorem}
Theorem \ref{t:1} provides the geometric interpretation we seek
for. The constant $3/16$ (absent in the Riemannian counterpart of
\eqref{e:11}) is due to the nonvanishing of $T_\nabla (X,Y)$ for
$X, Y \in H(M)$ i.e. of the torsion component proportional to the
Levi form. If in turn $\sigma \subset T_x (M)$ is a $2$-plane
tangent to $u = T_x$ and $\{ T_x , v \} \subset \sigma$ is a
$g_{\theta , x}$-orthonormal basis of $\sigma$ then we shall show
that
\begin{equation}
\label{e:12} L(\beta_r ) = 2\pi r - \frac{\pi r^3}{12} \{ 16
K_\theta (\sigma ) + \frac{3}{2}\; A_x (v,v)^2 +
\end{equation}
\[+ 2 \Omega_x (\tau_x v , v) - \| \tau_x v \|^2
\} + O(r^4 )\] where $A(X,Y) = g_\theta (X, \tau Y)$ and $\Omega =
- d \theta$. So an interpretation similar to that in Theorem
\ref{t:1} is not available unless $(M, \theta )$ is a Sasakian
manifold. Indeed if this is the case ($\tau = 0$) then we obtain
$K_\theta (\sigma ) = \lim_{r \to 0} (3/(4\pi r^3 )) \{ 2\pi r -
L(\beta_r )\}$.

\vskip 0.5cm

{\em Proof of Theorem} \ref{t:1}. Let $\sigma \subset H(M)_x$ be a
holomorphic $2$-plane and $v = J_x u$ where $u \in \sigma$, $\| u
\| = 1$. Recall (as a consequence of the purity axioms, cf. also
Chapter I in \cite{kn:DrTo}) that
\begin{equation}
T_\nabla (X,Y) = - \Omega (X,Y)T, \;\;\; X,Y \in H(M) \; ,
\label{e:13}
\end{equation}
hence $f^{\prime\prime\prime}(0) = 0$. On the other hand
\[ (\nabla_X T_\nabla )(Y,Z) = - (\nabla_X \Omega )(Y,Z) T = 0 \]
for any $X \in T(M)$, $Y,Z \in H(M)$, hence $P_s (u,v) = 0$. Also
(again by \eqref{e:13}) $\| T_{\nabla , x}(u,v)\| = 1$ and
\[ \langle T_{\nabla , x} (T_{\nabla , x} (u,v) ,
w(s)) \, , \, w(s + \frac{\pi}{2} )\rangle = \langle \tau_x (w(s))
\, , \, w(s + \frac{\pi}{2}) \rangle =
\]
\[ = A_x (w(s), w(s+ \frac{\pi}{2}) \rangle = - g(s) \]
where $g(s) = (\sin 2s ) A_x (u,u) - (\cos 2s ) A_x (u,v)$,
because of $A_x (v,v) = - A_x (u,u)$ (itself a consequence of
$\tau \circ J = - J \circ \tau$). It follows that \[ f^{(4)}(0) =
6 - 32 H_\theta (\sigma ) + 8 g(s) \; . \] Summing up $f(r) =
\sum_{j=0}^4 \frac{f^{(j)}(0)}{j!} r^j + O(r^5 ) = r^2 (1 - \delta
)$, where $\delta = (r^2 /12) \{ 16 H_\theta (\sigma ) - 3 - 4
g(s)\} + O(r^3 )$, hence
\[ \| X_s \|_{\gamma_{w(s)}(r)} = r \sqrt{1 - \delta}
= r(1 - \frac{\delta}{2} + O(\delta^2 )) = \]
\[ = r - \frac{r^3}{24} \{ 16 H_\theta (\sigma ) - 3
- 4 g(s) \} + O(r^4 ) \; .\] Finally, by integration we obtain (as
$\int_0^{2\pi} g(s) \, d s = 0$) the identity \eqref{e:3} and the
proof of Theorem \ref{t:1} is complete. For $2$-planes tangent to
$u = T_x$ we have
\[ f^{\prime\prime\prime} (0) = 6 A_x (v , w(s + \frac{\pi}{2})) \;.
\]
and
\[ \langle T_{\nabla , x}(T_{\nabla , x}(u,v), w(s)), w(s +
\frac{\pi}{2}) \rangle = \| \tau_x v \|^2 \cos^2 s + \Omega_x
(\tau_x v , v) \sin^2 s \; . \] Also $(\nabla_X T_\nabla )(T, Y) =
(\nabla_X \tau ) Y$ implies
\[ \langle P_{s,x} (u,v), w(s+ \frac{\pi}{2}) \rangle =
(\nabla_{\dot{\gamma}_{w(s)}} A)_x (v,v) \cos s \] hence
\[
f^{(4)}(0) = 6 \| \tau_x v \|^2 - 32 K_\theta (\sigma ) +
16(\nabla_{\dot{\gamma}_{w(s)}} A)_x (v,v) \cos s - \]
\[  - 8
\{ \| \tau_x v \|^2 \cos^2 s + \Omega_x (\tau_x v , v) \sin^2 s \}
\; . \] Similar to the above we set
\[ \delta = - r A_x (v,v) \cos s + \frac{r^2}{12} \{
16 K_\theta (\sigma ) + (4 \cos^2 s - 3) \| \tau_x v \|^2 -
\]
\[ - 8(\nabla_{\dot{\gamma}_{w(s)}} A)_x (v,v) \cos s + 4 \Omega_x
(\tau_x v , v) \sin^2 s \} + O(r^3 ) \] hence
\[ \| X_s \|_{\gamma_{w(s)}(r)} = r \sqrt{1 - \delta} = r(1 -
\frac{\delta}{2} - \frac{\delta^2}{8} + O(\delta^3 )) \] and
integration over $0 \leq s \leq 2 \pi$ leads to \eqref{e:12}.
\section{A Schur-like result} The scope of this section is to
establish the following
\begin{theorem} Let $M$ be a connected strictly pseudoconvex CR manifold of CR dimension
$n \geq 2$ and $\theta$ a contact form on $M$ with $L_\theta$
positive definite. Let $S(X,Y) = (\nabla_X \tau )Y - (\nabla_Y
\tau ) X$. Assume that the pseudohermitian sectional curvature is
a point function only i.e. $K_\theta = f \circ \pi$ for some
$C^\infty$ function $f : M \to {\mathbb R}$. If $S = 0$ then
$\nabla f = \theta (\nabla f) T$. Moreover if $(M , \theta )$ is a
Sasakian manifold $(\tau = 0)$ then $f = 0$; consequently $R = 0$
and $(M , \theta )$ is a Sasakian space form of sectional
curvature $c = -3$. \label{t:2}
\end{theorem}
Here $\nabla f$ is the ordinary gradient of $f$ with respect to
the Webster metric i.e. $g_\theta (\nabla f , X) = X(f)$ for any
$X \in T(M)$. As a byproduct of Theorem \ref{t:2} there are no
``pseudohermitian space forms'' except for those with $K_\theta =
0$. Moreover (as argued in \cite{kn:BaDr}) these aren't
Tanaka-Webster flat unless $\tau = 0$. So the term {\em
pseudohermitian space form} should be reserved for pseudohermitian
manifolds $(M , \theta )$ such that the sectional curvature
\eqref{e:1} (rather than \eqref{e:2}) is constant and then
examples abound. For instance (cf. \cite{kn:Web1} or Section 1.5
in \cite{kn:DrTo}) if $[g_{\alpha\overline{\beta}}] \in {\rm
GL}(n, {\mathbb C})$ is a Hermitian matrix and $c \in (0, + \infty
)$ then let $Q_\pm (c) \subset {\mathbb C}^{n+1}$ be the real
hypersurface defined by $r_\pm (z,w) \equiv
g_{\alpha\overline{\beta}} z^\alpha \overline{z}^\beta \pm (w
\overline{w} - c) = 0$, where $(z^1 , \cdots , z^n , w )$ are the
natural complex coordinates on ${\mathbb C}^{n+1}$. Then $Q_\pm
(c)$ is a nondegenerate CR manifold and the contact form
$\theta_\pm = i g_{\alpha\overline{\beta}} (z^\alpha d
\overline{z}^\beta - \overline{z}^\beta d z^\alpha ) \pm i (w d
\overline{w} - \overline{w} d w )$ has constant sectional
curvature ${\rm Sect}(Q_\pm (c), \theta_\pm ) = \pm 1/(2c)$. To
prove Theorem \ref{t:2} let us set
\[ R_1 (X,Y,Z,W) = g_\theta (X,Z) g_\theta (Y,W) - g_\theta (Y,Z)
g_\theta (X,W) \] for any $X,Y,Z,W \in T(M)$. If $L := R - 4 f
R_1$ then (by hypothesis)
\[ L (X,Y,X,Y) = 0, \;\;\; X,Y \in T(M) \; . \]
Thus (by a result in \cite{kn:BaDr}, Appendix A)
\begin{equation}
\label{e:14} R(X,Y,Z,W) = 4 f R_1 (X,Y,Z,W) + \Omega (Y,W) A(X,Z)
-
\end{equation}
\[ - \Omega (Y,Z) A(X, W)  + \Omega (X,Z) A(Y,W) - \Omega (X,W) A(Y,Z)
+ \]
\[+ g_\theta (S(Z_H , W_H ) \, , \, (\theta \wedge
I)(X,Y)) - g_\theta (S(X_H , Y_H ) \, , \, (\theta \wedge I)(Z,W))
\]
where $X_H = \pi_H X$ and $\pi_H : T(M) \to H(M)$ is the
projection associated to the decomposition $T(M) = H(M) \oplus
{\mathbb R}$. Also $I$ is the identical transformation and
$(\theta \wedge I)(X,Y) = \frac{1}{2} \{ \theta (X) Y - \theta (Y)
X \}$. Note that $\nabla g_\theta = 0$ yields $\nabla R_1 = 0$
hence (by computing the covariant derivative of \eqref{e:14} and
using $\nabla \Omega = 0$)
\begin{equation}
\label{e:15} (\nabla_U R)(X,Y,Z,W) = U(f) R_1 (X,Y,Z,W) +
\end{equation}
\[ + \Omega (Y,W) (\nabla_U A)(X,Z) - \Omega
(Y,Z) (\nabla_U A)(X,W) + \]
\[ + \Omega (X,Z)
(\nabla_U A)(Y,W) - \Omega (X,W) (\nabla_U A)(Y,Z) +  \]
\[ +
g_\theta ((\nabla_U S)(Z_H , W_H ) \, , \, (\theta \wedge I)(X,Y))
-\] \[-  g_\theta ((\nabla_U S)(X_H , Y_H ) \, , \, (\theta \wedge
I)(Z,W))\] for any $X,Y,Z,W,U \in T(M)$. Let us take the cyclic
sum over $(U,Z,W)$ and use the second Bianchi identity (cf.
Theorem 5.3 in \cite{kn:KoNo}, Vol. I, p. 135)
\[ \sum_{UZW} (\nabla_U R)(X,Y,Z,W) = - \sum_{UZW} g_\theta
(R(T_\nabla (U,Z), W)Y, X) \] so that to obtain
\begin{equation}
\label{e:16} - \sum_{UZW} R(T_\nabla (U,Z), W)Y =  U(f) \{g_\theta
(Y,W)Z - g_\theta (Y,Z) W \} + \end{equation}
\[  +
Z(f) \{ g_\theta (Y,U) W - g_\theta (Y,W) U \}  +\] \[ + W(f) \{
g_\theta (Y,Z) U - g_\theta (Y,U) Z \} + \]
\[ + \Omega (Y,W) S(U,Z) +
\Omega (Y, U) S(Z,W) + \Omega (Y,Z) S(W,U) +
\]
\[ + g_\theta (Y\, , \, S(U,W) J Z + S(Z,U) J W
+ S(W, Z) J U) - \]
\[  - g_\theta ((\nabla_U S)(\pi_H \, \cdot
\, , \, Y_H ) \, , \, (\theta \wedge I)(Z,W))^\sharp - \]
\[-
g_\theta ((\nabla_Z S)(\pi_H \, \cdot \, , \, Y_H ) \, , \,
(\theta \wedge I) (W,U))^\sharp - \]
\[ - g_\theta ((\nabla_W S)(\pi_H \, \cdot , Y_H ) \, , \, (\theta \wedge I)(U,Z))^\sharp - \]
\[- \frac{1}{2}\; \theta (Y) \{ (\nabla_U S)(Z_H ,
W_H ) + (\nabla_Z S)(W_H , U_H )  + (\nabla_W S)(U_H , Z_H ) \} +
\]
\[ + \frac{1}{2}\; g_\theta (Y \, , \,
(\nabla_U S)(Z_H , W_H ) + (\nabla_Z S)(W_H , U_H )  +(\nabla_W
S)(U_H , Z_H ) ) T \] for any $Y,Z,W,U \in T(M)$. Here $\sharp$
denotes raising of indices with respect to $g_\theta$ i.e.
$g_\theta (\omega^\sharp , X) = \omega (X)$ for any $\omega \in
T^* (M)$ and any $X \in T(M)$. In particular for $Y,Z,W,U \in
H(M)$ the left hand member of \eqref{e:16} becomes (by
\eqref{e:13}) $\sum_{UZW} \Omega (U,Z) R(T, W)Y$. To compute terms
of the form $R(T, Y)Z$ we need to recall the identity (cf. Section
1.4.2 in \cite{kn:DrTo})
\begin{equation}
\label{e:17} g_\theta (R(X,Y)Z , W) = g_\theta (R(W,Z)Y , X) -
g_\theta ((L X \wedge L Y)Z , W) + \end{equation}
\[ + g_\theta
((L W \wedge L Z)Y , X) +  g_\theta (S(X,Y) , Z) \theta (W) -
g_\theta (S(W,Z) , Y) \theta (X) - \]
\[  -
\theta (Z) g_\theta (S(X,Y) , W) + \theta (Y) g_\theta (S(W,Z) ,
X) + \]
\[  + 2 g_\theta ((\theta \wedge {\mathcal O})(X,Y) ,
Z) \theta (W) - 2 g_\theta ((\theta \wedge {\mathcal O})(W,Z) , Y)
\theta (X) -
\]
\[ - 2 \theta (Z) g_\theta ((\theta \wedge
{\mathcal O})(X,Y) , W) + 2 \theta (Y) g_\theta ((\theta \wedge
{\mathcal O})(W, Z), X)
\]
for any $X,Y,Z,W \in T(M)$. Here
\[ L = \tau + J, \;\;\; {\mathcal O} = \tau^2 + 2 J \tau - I \; . \]
Also $(X \wedge Y) Z = g_\theta (Z, X) Y - g_\theta (Z, Y) X$. The
lack of symmetry of $R(X,Y,Z,W)$ in the pairs $(X,Y)$ and $(Z,W)$
(in contrast with the case of Riemannian curvature, cf.
Proposition 2.1 in \cite{kn:KoNo}, p. 201) is the consequence of
the presence of torsion terms in the first Bianchi identity. Let
us set $X = T$ and $Y,Z,W \in H(M)$ in \eqref{e:17}. We obtain (as
$L T = 0$)
\begin{equation}
g_\theta (R(T, Y)Z , W) = g_\theta (Y, S(Z,W)) \; . \label{e:18}
\end{equation}
Next for any vector field $Z \in H(M)$ we may choose $Y \in H(M)$
such that $g_\theta (Y,Z) = 0$ and $\| Y \| = 1$. Also let $U = Y$
and $W = J Y$. Then \eqref{e:16} becomes (by \eqref{e:18})
\[ g_\theta (Z, S(Y, J Y)) = Z(f) - g_\theta (J Y, S(Y,Z)) +
g_\theta (Y, S(JY,Z)) \] or $Z(f) = 2 g_\theta (S(Y, J Y) , Z)$
yielding the first statement in Theorem \ref{t:2}. Similarly we
may use \eqref{e:16} for $Z = T$ and $W \in H(M)$ chosen such that
$\| W \| = 1$ together with $U = Y$ and $Y = J W$ so that to
obtain (when $\tau = 0$) $T(f) W - W(f) T = 0$. As $M$ is
connected $f$ is constant and then by Theorem 5 in \cite{kn:BaDr}
it follows that $R = 0$ [and in particular $M$ is a spherical CR
manifold i.e. the Chern-Moser tensor vanishes identically (cf.
\cite{kn:Kam}, p. 187)]. If this is the case then (by Proposition
4 in \cite{kn:BaDr}) $(M, (J, \, - T, -\theta , g_\theta ))$ is a
Sasakian space form of (constant) $\varphi$-sectional curvature $c
= - 3$. Finally, if $M$ is compact let $(\rho , {\rm dev}) : ({\rm
Aut}_{\rm CR}(\tilde{M}), \tilde{M}) \to ({\rm PU}(n+1,1),
S^{2n+1})$ be the developing pair for $M$ as a spherical CR
manifold (cf. \cite{kn:Kam}, p. 195) where $\tilde{M}$ is the
universal covering space of $M$. Then (cf. \cite{kn:Kam}, p. 205)
${\rm dev} : \tilde{M} \to S^{2n+1} \setminus \{ \infty \} \approx
{\mathbb H}_n$ is an isometry (where ${\mathbb H}_n$ is thought of
as carrying the left invariant Webster metric associated to the
contact form $\theta_0 = dt + i \sum_{j=1}^n \left( z^j d
\overline{z}^{j} - \overline{z}^{j} d z^j \right)$) thus proving
Corollary \ref{c:1} (cf. also Theorem 6.1 in \cite{kn:Kam}, p.
207).
\section{Pseudohermitian manifolds of maximal $\dim_{\mathbb R}{\mathfrak i}(M, \theta )$}
\subsection{Infinitesimal pseudohermitian transformations} The
purpose of this section is to prove the following
\begin{theorem} Let $(M , \theta )$ be a connected pseudohermitian manifold
of CR dimension $n$ with $L_\theta$ positive definite. Then {\rm
a)} $\dim_{\mathbb R} {\mathfrak i}(M , \theta ) \leq (n+1)^2$.
{\rm b)} If $\dim_{\mathbb R} {\mathfrak i}(M , \theta ) =
(n+1)^2$ then $(M , \theta )$ has constant pseudohermitian
sectional curvature $H_\theta (\sigma )$.  \label{t:3}
\end{theorem}
\par
Let $(M , \theta )$ a $(2n+1)$-dimensional pseudohermitian
manifold of CR dimension $n$. A {\em CR isomorphism} is a
$C^\infty$ diffeomorphism $f : M \to M$ and a CR map i.e. $(d_x f)
T_{1,0}(M)_x = T_{1,0}(M)_{f(x)}$ for any $x \in M$. A CR
isomorphism $f : M \to M$ is a {\em pseudohermitian
transformation} if $f^* \theta = \theta$. Let ${\rm Psh}(M, \theta
)$ be the group of all pseudohermitian transformations. By a
result of S. M. Webster, \cite{kn:Web1}, i) ${\rm Psh}(M , \theta
)$ is a Lie group of dimension $\leq (n+1)^2$ with isotropy groups
of dimension $\leq n^2$. Moreover ii) if $M$ is strictly
pseudoconvex then the isotropy groups are compact and if $M$ is
compact then ${\rm Psh}(M, \theta )$ is compact. The statement (i)
also follows from part (a) in Theorem \ref{t:3}. Indeed each
$1$-parameter subgroup of ${\rm Psh}(M, \theta )$ induces an
infinitesimal pseudohermitian transformation which is complete and
conversely, so that the Lie algebra of ${\rm Psh}(M, \theta )$ is
isomorphic to the Lie subalgebra of ${\mathfrak i}(M , \theta )$
consisting of all complete infinitesimal pseudohermitian
transformations. In particular, if $\dim_{\mathbb R} {\mathfrak
i}(M , \theta ) = (n+1)^2$ then $\dim {\rm Psh}(M, \theta ) =
(n+1)^2$ hence one may apply the classification (up to contact
homothetites) result Theorem 4.10 in \cite{kn:Mus}, p. 236.
Another proof of S. M. Webster's result (i)-(ii) above was given
by E. Musso, cf. {\em op. cit.}, p. 225.
\par
Let ${\rm GL}(m, {\mathbb R}) \to L(M)
\stackrel{\Pi}{\longrightarrow} M$ be the principal bundle of all
linear frames tangent to $M$, where $m = 2n+1$. Any diffeomorphism
$f : M \to M$ induces in a natural manner an automorphism
$\tilde{f}$ of ${\rm GL}(m, {\mathbb R}) \to L(M) \to M$ (cf. e.g.
\cite{kn:KoNo}, Vol. I, p. 226). Assume from now on that $M$ is
strictly pseudoconvex and $\theta$ is chosen such that $L_\theta$
is positive definite. Let $U(M , \theta )_x$ consist of all linear
frames $b \in L(M)_x$ such that
\[ b(e_0 ) = T_x \, , \;\;\; b(e_\alpha ) \in H(M)_x \, , \;\;\;
b(e_{\alpha + n}) = J_x b(e_\alpha ), \;\;\; 1 \leq \alpha \leq n
\; , \]
\[ g_{\theta , x}(b(e_i ) , b(e_j )) = \delta_{ij} \, , \;\;\; 0
\leq i,j \leq 2n \; . \] This construction gives rise to a
principal subbundle ${\rm U}(n) \to U(M, \theta )$ $\to M$ of
$L(M)$. By a result of S. Nishikawa et al. (cf. Proposition 10 in
\cite{kn:DrNi}, p. 1065) a diffeomorphism $f : M \to M$ is a
pseudohermitian transformation if and only if $\tilde{f}(U(M,
\theta )) = U(M, \theta )$. Also for any fibre-preserving
diffeomorphism $F : U(M, \theta ) \to U(M, \theta )$ leaving
invariant the canonical form $\nu$ ($\nu_b = b^{-1} \circ (d_b \Pi
)$, $b \in U(M, \theta )$) there is $f \in {\rm Psh}(M , \theta )$
such that $F = \tilde{f}$.
\par
A tangent vector field $X$ on $M$ is an {\em infinitesimal
pseudohermitian transformation} of $(M , \theta )$ if the local
$1$-parameter group of local transformations induced by $X$
consists of local pseudohermitian transformations of $(M , \theta
)$. Let $X$ be a vector field on $M$ and $\{ \varphi_t \}_{|t| <
\epsilon}$ the local $1$-parameter group of local transformations
induced by $X$. Let $\tilde{X}$ be the {\em natural lift} of $X$
to $L(M)$ (cf. \cite{kn:KoNo}, Vol. I, p. 229-230) i.e. the vector
field $\tilde{X}$ on $L(M)$ induced by the local $1$-parameter
group $\{ \tilde{\varphi}_t \}_{|t| < \epsilon}$ of local
transformations of $L(M)$. Let ${\mathfrak i}(M , \theta )$ denote
the set of all infinitesimal pseudohermitian transformations of
$(M , \theta )$. By Proposition 11 in \cite{kn:DrNi}, p. 1066, the
following statements are equivalent 1) $X \in {\mathfrak i}(M ,
\theta )$, 2) $\tilde{X}_b \in T_b (U(M, \theta ))$ for any $b \in
U(M, \theta )$, 3) ${\mathcal L}_X \theta = 0$ and ${\mathcal L}_X
\theta^\alpha = f^\alpha_\beta \, \theta^\beta$ for any local
frame $\{ \theta^\alpha : 1 \leq \alpha \leq n \}$ of
$T_{1,0}(M)^*$ defined on the open subset $U \subseteq M$ and some
$C^\infty$ functions $f^\alpha_\beta : U \to {\mathbb R}$. Here
${\mathcal L}_X$ denotes the Lie derivative. As a corollary of
${\mathcal L}_{[X,Y]} = [{\mathcal L}_X , {\mathcal L}_Y ]$ and
the previous characterization of ${\mathfrak i}(M , \theta )$ it
follows that ${\mathfrak i}(M , \theta )$ is a Lie algebra.
\subsection{Affine transformations} To prove the first statement
in Theorem \ref{t:3} it suffices to show that, for a fixed linear
frame $b \in U(M, \theta )$, the linear map \[ \Phi_b : {\mathfrak
i}(M , \theta ) \to T_b (U(M, \theta )), \;\;\; \Phi_b (X) =
\tilde{X}_b \, , \;\;\; X \in {\mathfrak i}(M , \theta ) \; , \]
is injective. Indeed \[ \dim_{\bf R} T_b (U(M, \theta )) = \dim M
+ \dim {\rm U}(n) = (n+1)^2 \; . \] Let $\nabla$ be the
Tanaka-Webster connection of $(M , \theta )$. An {\em affine
transformation} of $(M , \nabla )$ is a diffeomorphism $f : M \to
M$ such that $\nabla_{\dot{\gamma}_f} X = 0$ along $\gamma_f := f
\circ \gamma$, for any tangent vector field $X$ along $\gamma$
such that $(\nabla_{\dot{\gamma}} X )_{\gamma (t)} = 0$, and for
any curve $\gamma$ in $M$. Let ${\mathfrak U}(M, \nabla )$ be the
group of all affine transformations of $(M , \nabla )$. If $f : M
\to M$ is a diffeomorphism and $X$ is a vector field on $M$ we set
$(f_* X)_y = (d_{f^{-1}(y)} f)X_{f^{-1}(y)}$ for any $y \in M$. By
a result of J. Masamune et al. (cf. the proof of Lemma 1 in
\cite{kn:DrMa}, p. 357) if we set
\[ \nabla^f_X Y := (f_* )^{-1} \nabla_{f_* X} f_* Y \]
and $f$ is a pseudohermitian transformation then $\nabla^f =
\nabla$. Therefore we may apply Proposition 1.4 in \cite{kn:KoNo},
Vol. I, p. 228, to conclude that $f$ is an affine transformation,
hence ${\rm Psh}(M , \theta )$ is a subgroup of ${\mathfrak U}(M,
\nabla )$.
\par
A tangent vector field $X$ on $M$ is an {\em infinitesimal affine
transformation} of $(M , \nabla )$ if the local $1$-parameter
group induced  by $X$ consists of local affine transformations of
$(M , \nabla )$. Let ${\mathfrak a}(M, \nabla )$ be the Lie
algebra of all affine transformations of $(M , \nabla )$.
\par
Let $\omega \in \Gamma^\infty (T^* (L(M)) \otimes {\mathfrak
g}{\mathfrak l}(m, {\mathbb R}))$ be the connection $1$-form
associated to the Tanaka-Webster connection $\nabla$ and let us
denote by ${\mathfrak a}(\omega )$ the Lie algebra of all tangent
vector fields ${\mathcal X}$ on $L(M)$ such that 1) $(d_u R_a )
{\mathcal X}_u = {\mathcal X}_{ua}$, $u \in L(M)$, $a \in {\rm
GL}(m, {\mathbb R})$, 2) ${\mathcal L}_{\mathcal X} \nu = 0$, and
3) ${\mathcal L}_{\mathcal X} \omega = 0$. Here $\nu_b = b^{-1}
\circ (d_b \Pi )$ for any $b \in L(M)$. It is a well known fact
(cf. e.g. \cite{kn:KoNo}, Vol. I, p. 232) of general connection
theory that the map $X \mapsto \tilde{X}$ gives a Lie algebra
isomorphism ${\mathfrak a}(M, \nabla ) \approx {\mathfrak
a}(\omega )$.
\par
Now we may prove Theorem \ref{t:3}. To this end let $X \in {\rm
Ker}(\Phi_u )$. Then $X \in {\mathfrak i}(M, \theta ) \subset
{\mathfrak a}(M, \nabla )$ hence $\tilde{X} \in {\mathfrak
a}(\omega )$ and $\tilde{X}_u = 0$ hence one may apply the lemma
in \cite{kn:KoNo}, Vol. I, p.  232 (in the proof of Theorem 2.3,
cf. {\em op. cit.}) to conclude that $\tilde{X} = 0$ identically
on $L(M)$. Yet (by Proposition 2.1 in \cite{kn:KoNo}, Vol. I, p.
229) $\tilde{X}$ is $\Pi$-related to $X$ so $X = 0$ everywhere on
$M$.
\par
To prove the second statement in Theorem \ref{t:3} let $\sigma \in
G_2 (H(M))_x$ and let $b \in U(M , \theta )$ such that $\Pi (b ) =
x$. Let $u \in \sigma$ such that $\| u \| = 1$ and $\xi \in
{\mathbb C}^n$ given by $\xi = b^{-1}(u)$. Here ${\mathbb C}^n
\approx {\mathbb R}^{2n} \times \{ 0 \} \subset {\mathbb R}^m$.
Let $B(\xi )$ and $B(J_0 \xi )$ be the standard horizontal vector
fields associated (in the sense of \cite{kn:KoNo}, Vol. I, p. 119)
to $\xi$ and $J_0 \xi$, where $J_0$ is the standard complex
structure on ${\mathbb C}^n$. Let $\Omega = D \omega$ be the
curvature $2$-form of the Tanaka-Webster connection. Then, again
by a general fact within connection theory (cf. \cite{kn:KoNo},
Vol. I, p. 133)
\[ H_\theta (\sigma ) = g_{\theta , x}( R_x (u, J_x u)J_x u \, , \, u )  = 2
\left( \Omega (B(\xi ), B(J_0 \xi ))_b \cdot J_0 \xi \, , \, \xi
\right)\] where $(\; , \; )$ is the Euclidean inner product on
${\mathbb R}^m$ and $A \cdot \xi$ is the matrix product ($A \in
{\mathfrak g}{\mathfrak l}(m , {\mathbb R}) \approx {\mathbb
R}^{m^2}$).
\par
We wish to show that $H_\theta (\sigma )$ is a point function
only. To this end let $\sigma^\prime \in G_2 (H(M))_x$ be another
holomorphic frame tangent at $x \in M$ and $v \in \sigma^\prime$
such that $\| v \| = 1$. We set $\eta = b^{-1} (v) \in {\mathbb
C}^n$. There is $g \in {\rm U}(n)$ such that $\eta = g \xi$. Then
(by Proposition 2.2 in \cite{kn:KoNo}, Vol. I, p. 119)
\[ \Omega (B(\eta ), B(J_0 \eta ))_b = \Omega (B(g \xi ), B(J_0 g \xi
))_b = \]
\[= \Omega_b ((d_{bg} R_{g^{-1}} )B(\xi )_{b g} \, ,
\, (d_{bg} R_{g^{-1}} )B(J_0 \xi )_{b g}) = \]
\[= {\rm ad}(g)
\Omega_{b g} (B(\xi )_{b g} \, , \, B(J_0 \xi )_{b g} ) = g \cdot
\Omega_{b g} (B(\xi )_{b g} \, , \, B(J_0 \xi )_{b g}) \cdot
g^{-1}\] where $R_g : U(M, \theta ) \to U(M, \theta )$ is the
right translation by $g$ and ${\rm ad}$ denotes the adjoint
representation of ${\rm GL}(m, {\mathbb R})$ in its Lie algebra.
Moreover
\[
H_\theta (\sigma^\prime ) = 2 (\Omega (B(\eta ) , B(J_0 \eta ))_b
\cdot J_0 \eta \, , \, \eta ) = \]
\[ = 2 ((g \cdot \Omega_{b
g} (B(\xi )_{b g} \, , \, B(J_0 \xi )_{b g} ) \cdot g^{-1} ) \cdot
g J_0 \xi \, , \, g \xi ) = \]
\[ = 2 (\Omega_{b g} (B(\xi )_{b
g} \, , \, B(J_0 \xi )_{b g}) \cdot J_0 \xi \, , \, \xi )
\] and it remains to be shown that the function $\left.
F_x : = F \right|_{U(M, \theta )_x}$ is constant, where \[ F :
U(M, \theta ) \to {\mathbb R}, \;\,\; F(b) = \Omega (B(\xi ),
B(J_0 \xi ))_b \, , \;\;\; b \in U(M, \theta ) \; .\] Let $X \in
{\mathfrak i}(M, \theta ) \subset {\mathfrak a}(M, \nabla )$. Then
(by Proposition 2.2 in \cite{kn:KoNo}, Vol. I, p. 230)
\[
\tilde{X}(F) = \tilde{X}(\Omega (B(\xi ), B(J_0 \xi ))) =
({\mathcal L}_{\tilde{X}} \Omega )(B(\xi ), B(J_0 \xi )) +
\]
\[  + \Omega ([\tilde{X}, B(\xi )] , B(J_0 \xi )) + \Omega
(B(\xi ), [\tilde{X}, B(J_0 \xi )]) = 0 \; . \] This simple fact
has two consequences. First, let $V$ be an arbitrary tangent
vector on $U(M, \theta )_x$ i.e. \[ V \in T_b (U(M, \theta )_x ) =
{\rm Ker}(d_b \Pi ) \subset T_b (U(M, \theta )) \] for some $b \in
U(M, \theta )$ with $\Pi (b) = x$. As $\Phi_b$ is assumed to be
on-to there is $X \in {\mathfrak i}(M, \theta )$ such that
$\tilde{X}_b = V$ hence $V(F) = \tilde{X}(F)_b = 0$. As $U(n)$ is
connected (and $U(M, \theta )_x \approx U(n)$, a diffeomorphism)
it follows that $F$ is constant i.e. $\Omega_{b g}(B(\xi )_{b g} ,
B(J_0 \xi )_{b g}) = \Omega_b (B(\xi )_b , B(J_0 \xi )_b )$ hence
there is a smooth function $f : M \to {\mathbb R}$ such that
$H_\theta = f \circ \pi$. At this point we may apply Theorem
\ref{t:5} in Appendix A provided that $n \geq 3$ and $\tau = 0$.
However one may prove Theorem \ref{t:3} in full generality as
follows. Let $W \in T_b (U(M, \theta ))$ be an arbitrary tangent
vector and $Y \in {\mathfrak i}(M, \theta )$ such that
$\tilde{Y}_b = W$. Then $W(F) = \tilde{Y}(F)_b = 0$ hence for any
fixed $\xi \in {\mathbb C}^n$ the function $(\Omega (B(\xi ),
B(J_0 \xi ) )\cdot J_0 \xi \, , \, \xi )$ is constant in a
neighborhood of $b$, so that $f$ follows to be locally constant,
and then constant on $M$.
\section{Pseudohermitian space forms}
A pseudohermitian manifold $(M , \theta )$ with $H_\theta (\sigma
) =$ const. is said to be a {\em pseudohermitian space form}.
Similarly to Theorem 5 in \cite{kn:BaDr} (giving the precise form
of the curvature tensor field $R$ of $(M , \theta )$ when
$K_\theta (\sigma ) =$ const.) we establish
\begin{theorem} Let $(M , \theta )$ be a pseudohermitian manifold
of CR dimension $n$. If $H_\theta (\sigma ) = c$ {\rm (}with $c
\in {\mathbb R}${\rm )} for any $\sigma \in G_2 (H(M))$ then
\begin{equation}
\label{e:19} R(X,Y,Z,W) = c \{ 2 \Omega (X,Y) \Omega (Z,W) +
\end{equation}
\[ + g_\theta (X,Z) g_\theta (Y,W) - g_\theta
(X,W) g_\theta (Y,Z) +  \]
\[ + \Omega (X,Z)
\Omega (Y,W) - \Omega (X,W) \Omega (Y,Z) \} +  \]
\[ + g_\theta
(X,Z) A(Y, J W) - g_\theta (X,W) A(Y, J Z) + \]
\[ + g_\theta(Y,W) A(X, J Z) - g_\theta
(Y,Z) A(X, J W) + \]
\[ + \Omega (X,Z) A(Y,W) -
\Omega (X,W) A(Y,Z) + \]\[  + \Omega (Y,W) A(X,Z) - \Omega (Y,Z)
A(X,W)
\] for any $X,Y,Z,W \in H(M)$. In particular
\begin{equation}
{\rm Ric}(X,Y) = 2 c (n+1)g_\theta (X,Y) + 2(n-1)A(X,J Y)
\label{e:20}
\end{equation}
for any $X,Y \in H(M)$ hence each pseudohermitian space form $(M ,
\theta )$ is a pseudo-Einstein manifold of constant
pseudohermitian scalar curvature $\rho = 2cn(n+1)$. \label{t:4}
\end{theorem}
Here ${\rm Ric}(X,Y) = {\rm trace} \{ Z \mapsto R(Z,Y)X \}$. If
$\{ T_\alpha : 1 \leq \alpha \leq n \}$ is a local frame of
$T_{1,0}(M)$ on the open set $U \subseteq M$ then we set
$g_{\alpha\overline{\beta}} = L_\theta (T_\alpha ,
T_{\overline{\beta}})$ and $R_{\alpha\overline{\beta}} = {\rm
Ric}(T_\alpha , T_{\overline{\beta}})$. Then
$R_{\alpha\overline{\beta}}$ is the {\em pseudohermitian Ricci
tensor} and $\rho = g^{\alpha\overline{\beta}}
R_{\alpha\overline{\beta}}$ is the {\em pseudohermitian scalar
curvature}. Cf. J. M. Lee, \cite{kn:Lee} (or \cite{kn:DrTo},
Chapter 5) a contact form $\theta$ is {\em pseudo-Einstein} if
$R_{\alpha\overline{\beta}} = (\rho /n)
g_{\alpha\overline{\beta}}$. To prove Theorem \ref{t:4} we
consider the $4$-tensor field
\begin{equation}
\label{e:21} R_0 (X,Y,Z,W) = \frac{1}{4} \{ g_\theta (X,Z)
g_\theta (Y,W) - g_\theta (X,W) g_\theta (Y,Z) +
\end{equation}
\[ + \Omega (X,Z) \Omega (Y,W) - \Omega (X,W)
\Omega (Y,Z) + 2 \Omega (X,Y) \Omega (Z,W) \}
\]
for any $X,Y,Z,W \in H(M)$ and set $L = R - 4 c R_0$. Then we
exploit the symmetries of $L$ to establish \eqref{e:19} (using the
algebraic machinery in the proof of Proposition 7.1 in
\cite{kn:KoNo}, Vol. II, p. 166). Details are given in Appendix A
where we also prove a Sasakian version of the complex Schur
theorem.
\par
E. Musso has classified (cf. \cite{kn:Mus}) up to contact
homotheties the $G$-homogeneous pseudo-Einstein manifolds $(M ,
\theta )$ with $L_\theta$ positive definite. The same problem when
$L_\theta$ is but nondegenerate is open. We recall that a
pseudohermitian manifold is $G$-{\em homogeneous} if there is a
closed subgroup $G \subset {\rm Psh}(M, \theta )$ such that $G$
acts transitively on $M$. Also a {\em contact homothety} among two
pseudohermitian manifolds $(M , \theta )$ and $(M^\prime ,
\theta^\prime )$ is a CR diffeomorphism $f : M \to M^\prime$ such
that $f^* \theta^\prime = r \, \theta$ for some $r \in (0, +
\infty )$.
\par
Let $(M , \theta )$ be a $G$-homogeneous pseudohermitian manifold
with $G$ connected and $L_\theta$ positive definite. As usual we
fix a point $x_0 \in M$ and let $H \subset G$ be the isotropy
subgroup at $x_0$ and $H \to G \to M = G/H$ the corresponding
principal bundle. Let $V$ be the left invariant vector field on
$G$ determined by $T_{x_0}$. Let $\mathfrak g$ and $\mathfrak h$
be the Lie algebras of $G$ and $H$, respectively. We consider a
reductive decomposition $\mathfrak{g} = \mathfrak{h} \oplus
\mathfrak{p}$ where $\mathfrak p$ is identified with $T_{x_0}
(M)$. Due to this identification one has a direct sum
decomposition ${\mathfrak p} = {\mathfrak m} \oplus {\mathfrak v}$
where the $H$-invariant subspaces $\mathfrak m$ and $\mathfrak v$
correspond to $H(M)_{x_0}$ and ${\mathbb R} T_{x_0}$,
respectively. Let $\eta$ be the left invariant differential
$1$-form on $G$ determined by
\[ ({\mathfrak h} \oplus {\mathfrak m}) \, \rfloor \, \eta = 0,
\;\;\; \eta (V) = 1 \; , \] and let us set $K = \{ a \in G : {\rm
ad}(a)^* \eta = \eta \}$. Finally let $K^\prime = K_0 H$, where
$K_0$ is the connected component of the identity in $K$, and $B =
G/K^\prime$. Then the natural projection $p : M \to B$ organizes
$M$ as a principal bundle (with $S^1$ or $\mathbb R$ as a
structure group) over $B$ (and the fibres of $p$ are maximal
integral curves of $T$). Combining Theorem \ref{t:4} above with
Theorems 4.5, 4.7 and 4.8 in \cite{kn:Mus}, p. 233-236, we may
conclude that Corollary \ref{c:2} holds.
\par
Let us briefly describe the pseudohermitian structures on the
model spaces in (i)-(iii) of Corollary \ref{c:2}. Under the
assumptions of Corollary \ref{c:2} it follows (by \eqref{e:20})
that $B$ is a simply-connected compact homogeneous
K\"ahler-Einstein manifold. Therefore, by a result in
\cite{kn:Mus}, p. 230-232, there is a principle $S^1$-bundle
$\pi_{(1)} : B_{(1)} \to B$ and a canonical contact form
$\theta_{(1)}$ on $B_{(1)}$ such that $(B_{(1)} , \theta_{(1)})$
is a pseudohermitian manifold with $L_{\theta_{(1)}}$ positive
definite and $\pi_{(1)}$ is a Riemannian submersion of $(B_{(1)} ,
g_{\theta_{(1)}})$ onto $B$. Moreover if $c_1 (M)$ is the integral
first Chern class of $T(M)$ then $c_1 (M) = k \, c_1 (B_{(1)})$
for some $k \in {\mathbb Z}$, $k > 0$. Let $\pi_{(k)} : B_{(k)}
\to B$ be the $k$th tensor power of $\pi_{(1)} : B_{(1)} \to B$.
Again by a result in \cite{kn:Mus}, p. 232, there is a unique
pseudohermitian structure $\theta_{(k)}$ on $B_{(k)}$ such that
$(B_{(k)} , \theta_{(k)})$ is a pseudohermitian manifold with
$L_{\theta_{(k)}}$ positive definite and $\pi_{(k)}$ is a
Riemannian submersion of $(B_{(k)} , g_{\theta_{(k)}})$ onto $(B,
\sqrt{k} g)$, where $g$ is the K\"ahler-Einstein metric of $B$.
Then $(B_{(k)} , \theta_{(k)})$ is referred to as the {\em
canonical pseudohermitian manifold of index} $k$ over $B$. The
contact form of the model space $B \times S^1$ in (ii) is given by
$\theta^\prime = a \, d \gamma + i(\overline{\partial} - \partial
) \log K(z, z)$ for some $a \in (0, + \infty )$, where $\gamma$ is
a local fibre coordinate (i.e. $\partial /\partial \gamma$ is
tangent to the $S^1$-action on $B \times S^1$) and $K(z, \zeta )$
is the Bergman kernel of $B$ (thought of as an affinely
homogeneous Siegel domain of the second kind, cf. Theorem 4.7 in
\cite{kn:Mus}, p. 235). Similarly $B \times {\mathbb R}$ is
endowed with the contact form $\theta^{\prime\prime} = a \, d t +
i(\overline{\partial} -
\partial ) \log K(z,z)$ for some $a \in (0, + \infty )$. As to the
model spaces in (iii), ${\mathbb C}^n \times S^1$ is endowed with
the contact form $\theta^\prime = a \, d \gamma + 2 \sum_{j=1}^n
y^j \, d x^j$ while ${\mathbb C}^n \times {\mathbb R}$ carries
$\theta^{\prime\prime} = a \, dt + 2 \sum_{j=1}^n y^j \, d x^j$.
\section{Pseudohermitian immersions}
Let $M$ and $M^\prime$ be two CR manifolds of CR dimensions $n$
and $n+k$ respectively, with $k \geq 1$. A {\em CR immersion} is a
$C^\infty$ immersion $f : M \to M^\prime$ and a CR map. Given
pseudohermitian structures $\theta$ and $\theta^\prime$ on $M$ and
$M^\prime$ respectively, a CR immersion is {\em
isopseudohermitian} if $f^* \theta^\prime = \theta$. Assume that
$M$ and $M^\prime$ are nondegenerate and let $T^\prime$ be the
characteristic direction of $d \theta^\prime$. A {\em
pseudohermitian immersion} is an isopseudohermitian CR immersion
such that $T^{\prime \bot} = 0$. If $V \in T(M^\prime )$ then
$V^\bot = {\rm nor}(V)$ and ${\rm nor} : T(M^\prime ) \to E(f)$ is
the projection associated to the decomposition $T(M^\prime ) =
[f_* T(M)] \oplus E(f)$ while $E(f) \to M$ denotes the normal
bundle of the given immersion. Here we assume that $(d_x f) T_x
(M)$ is nondegenerate in $(T_{f(x)}(M^\prime ), g_{\theta^\prime ,
f(x)})$ and then $E(f)_x$ is the $g_{\theta^\prime ,
f(x)}$-orthogonal complement of $(d_x f) T_x (M)$. A theory of
pseudohermitian immersions has been started by S. Dragomir,
\cite{kn:Dra}. Cf. also \cite{kn:BaDr2}. Assume from now on that
both $M$ and $M^\prime$ are strictly pseudoconvex and $\theta$,
$\theta^\prime$ are chosen such that $L_\theta$,
$L_{\theta^\prime}$ are positive definite. We shall need the
pseudohermitian analogs of the Gauss and Weingarten formulae
\begin{equation}
\nabla^\prime_{f_* X} f_* Y = f_* \nabla_X Y + \alpha (f)(X,Y) \;
, \label{e:Gauss}
\end{equation}
\begin{equation}
\nabla^\prime_{f_* X} \xi = - f_* a_\xi X + \nabla^\bot_X \xi \; ,
\label{e:Weingarten}
\end{equation}
for any $X,Y \in T(M)$ and any $\xi \in \Gamma^\infty (E(f))$. Cf.
(41)-(42) in \cite{kn:Dra}, p. 185. Here $\nabla^\prime$ is the
Tanaka-Webster connection of $(M^\prime , \theta^\prime )$ while
$\alpha (f)$ is a $E(f)$-valued $C^\infty (M)$-bilinear form,
$a_\xi$ is an endomorphism of $T(M)$, and $\nabla^\bot$ is a
connection in $E(f) \to M$ (the pseudohermitian analogs to the
second fundamental form, Weingarten operator and normal connection
of an isometric immersion). Let $R^\prime$ be the curvature tensor
field of $\nabla^\prime$. We recall (cf. (61) in \cite{kn:Dra}, p.
191)
\[ {\rm tan} \{ R^\prime (f_* X , f_* Y ) f_* Z \} = R(X,Y)Z +
a_{\alpha (f)(X,Z)} Y - a_{\alpha (f)(Y,Z)} X \] for any $X,Y,Z
\in T(M)$, where ${\rm tan} : T(M^\prime ) \to T(M)$ is the
natural projection. Let us take the inner product with $W \in
T(M)$ and use
\[ g_{\theta^\prime}(\alpha (f)(X,Y) , \xi ) = g_\theta (a_\xi X ,
Y) \] (cf. (50) in \cite{kn:Dra}, p. 188) so that to get
\begin{equation}
\label{e:curv} R^\prime (f_* W , f_* Z , f_* X , f_* Y) =
R(W,Z,X,Y) + \end{equation}
\[+ g_{\theta^\prime} (\alpha
(f)(X,Z) \, , \, \alpha (f)(Y,W))  - g_{\theta^\prime} (\alpha
(f)(Y,Z), \alpha (f)(X,W)) \; .
\]
\begin{lemma} For any $X,Y \in T(M)$
\begin{equation}
\alpha (f) (X, J Y) = J^\prime \alpha (f)(X,Y) \; ,
\label{e:fund1}
\end{equation}
\begin{equation}
\alpha (f)(J X , Y) = J^\prime \alpha (f)(X,Y) - \theta (X)
J^\prime Q Y
\label{e:fund2}
\end{equation}
where $J^\prime$ is the complex structure on $H(M^\prime )$ {\rm
(}extended to an endomorphism of $T(M^\prime )$ by requiring that
$J^\prime T^\prime = 0${\rm )} and $Q (X) = \alpha (f)(T, X)$.
Consequently
\begin{equation} \alpha (J X , J Y) = - \alpha
(f)(X,Y) + \theta (X) Q Y \; . \label{e:fund3}
\end{equation}
for any $X,Y \in T(M)$. \label{l:1}
\end{lemma}
\noindent Replacing $(W,Z,X,Y)$ by $(X, J X, X, J X)$ in
\eqref{e:curv} and using Lemma \ref{l:1} leads to the following
\begin{theorem} Let $f : M \to M^\prime$ be a pseudohermitian
immersion between two pseudohermitian manifolds $(M , \theta )$
and $(M^\prime , \theta^\prime )$ such that $L_\theta$ and
$L_{\theta^\prime}$ are positive definite. Then
\[
R^\prime (f_* X , J^\prime f_* X , f_* X , J^\prime f_* X) =
R(X,JX,X,JX) + \]
\[ + 2 \| \alpha (f)(X,X)\|^2 - 2 \theta (X)
g_{\theta^\prime}(\alpha (f)(X,X) \, , \, Q X) \] for any $X \in
T(M)$. In particular $H_\theta (\sigma ) \leq
H_{\theta^\prime}((d_x f) \sigma )$ for any $\sigma \in G_2
(H(M))_x$ and any $x \in M$. \label{t:6}
\end{theorem}
It remains that we prove Lemma \ref{l:1}. The identity
\eqref{e:fund1} is a consequence of $\nabla^\prime J^\prime = 0$
and the Gauss formula \eqref{e:Gauss}. Cf. also (43) in
\cite{kn:Dra}, p. 187. Moreover the identity
\[ T_{\nabla^\prime} = 2 (\theta^\prime \wedge \tau^\prime -
\Omega^\prime \otimes T^\prime ) \] (cf. e.g. \cite{kn:DrTo},
Chapter 1) leads to
\begin{equation} \alpha (f)(Y,X) = \alpha
(f)(X,Y) - 2 (\theta \wedge Q)(X,Y)
\label{e:symm}
\end{equation}
where $Q(X) = \alpha (f)(T, X)$ for any $X \in T(M)$. Finally
\eqref{e:fund1} and \eqref{e:symm} imply
\eqref{e:fund2}-\eqref{e:fund3}.
\par
The proof of Corollary \ref{c:3} is immediate. Corollary \ref{c:4}
follows from a result by L. Ornea \& M. Verbitsky, cf. Theorem 6.1
in \cite{kn:OrVe}, p. 141. Indeed let $(M , \theta )$ be a compact
Sasakian manifold and $V = M \times S^1$. Then $V$ is a Vaisman
manifold (cf. e.g. \cite{kn:DrOr} for the relevant notions)
admitting (cf. Theorem 5.1 in \cite{kn:OrVe}, p. 138) an immersion
$\phi : M \to H_\Lambda$ into a primary Hopf manifold $H_\Lambda =
({\mathbb C}^{n+1} \setminus \{ 0 \} )/\Gamma_\Lambda$ for some $n
\geq 1$ and some $\Lambda = (\lambda_1 , \cdots , \lambda_{n+1} )
\in {\mathbb C}^{n+1}$ such that $0 < |\lambda_{n+1} | \leq \cdots
\leq |\lambda_1 | < 1$. Here $\Gamma_\Lambda$ is the discrete
group of complex analytic transformations of ${\mathbb C}^{n+1}
\setminus \{ 0 \}$ generated by $(z_1 , \cdots , z_{n+1} ) \mapsto
(\lambda_1 z_1 , \cdots , \lambda_{n+1} z_{n+1} )$. See also
\cite{kn:KaOr}, p. 202. Moreover $\phi$ descends to a
pseudohermitian immersion $M \to (S^{2n+1} , \theta_A )$ with
$\lambda_j = e^{-a_j}$ hence (by Theorem \ref{t:6} above) the
upper bound on $H_\theta (\sigma )$ in Corollary \ref{c:4}.

\begin{proposition} Let $(M , \theta )$ be a pseudohermitian
manifold with $L_\theta$ positive definite. If $\hat{\theta} =
e^{2u} \theta$, $u \in C^\infty (M)$, then
\begin{equation}
\label{e:relsect}
e^{2u} H_{\hat{\theta}} (\sigma ) =
H_\theta (\sigma ) + 2 i u_0 - 2 u_\alpha u^\alpha - 2
(\nabla_{\overline{\beta}} u_\alpha ) \eta^\alpha
\eta^{\overline{\beta}}
\end{equation}
for any $\sigma \in G_2 (H(M))_x$ and $x \in M$, where $X = Z +
\overline{Z} \in \sigma$, $Z = \xi^\alpha T_\alpha$, and
$\eta^\alpha = \| \xi \|^{-1} \xi^\alpha$, $\| \xi \|^2 =
g_{\alpha\overline{\beta}} \xi^\alpha \xi^{\overline{\beta}}$.
Consequently the pseudohermitian sectional curvature is not a CR
invariant. In particular if $\hat{\theta} = (1/a) \theta$ {\rm
(}$a > 0${\rm )} then $H_{\hat{\theta}} (\sigma ) = a H_\theta
(\sigma )$. \label{p:1}
\end{proposition}
{\em Proof}. Let $\{ T_\alpha : 1 \leq \alpha \leq n \}$ be a
local frame of $T_{1,0}(M)$. Let $\hat{\nabla}$ be the
Tanaka-Webster connection of $(M , e^{2u} \theta )$ and
$\hat{\Gamma}^A_{BC}$ the connection coefficients with respect to
$\{ T_A : A \in \{ 0, 1, \cdots , n , \overline{1} , \cdots ,
\overline{n} \} \}$ (with the convention $T_0 = T$). We set
${{R_C}^D}_{AB} T_D  = R(T_A , T_B ) T_C$ so that
\[
{{R_\alpha}^\beta}_{\lambda \overline{\mu}} = T_\lambda
(\Gamma^\beta_{\overline{\mu}\alpha}) -
T_{\overline{\mu}}(\Gamma^\beta_{\lambda\alpha}) + 2 i
\Gamma_{0\alpha}^\beta g_{\lambda\overline{\mu}} + \]
\[ +
\Gamma^\gamma_{\overline{\mu}\alpha} \Gamma^\beta_{\lambda\gamma}
- \Gamma^\gamma_{\lambda\alpha}
\Gamma^\beta_{\overline{\mu}\gamma} +
\Gamma^\gamma_{\overline{\mu}\lambda} \Gamma^\beta_{\gamma\alpha}
- \Gamma^{\overline{\gamma}}_{\lambda\overline{\mu}}
\Gamma^\beta_{\overline{\gamma}\alpha} \; . \] The proof of
Proposition \ref{p:1} is to replace in
${{{\hat{R}}\,_\alpha}^\beta}_{\lambda\overline{\mu}}$ from the
identities
\begin{equation} \hat{\Gamma}^\alpha_{\gamma\beta} =
\Gamma^\alpha_{\gamma\beta} + 2(u_\gamma \delta^\alpha_\beta +
u_\beta \delta^\alpha_\gamma ) \; , \label{e:coeff1}
\end{equation}
\begin{equation}
\hat{\Gamma}^\alpha_{\overline{\gamma}\beta} =
\Gamma^\alpha_{\overline{\gamma}\beta} - 2 u^\alpha
g_{\beta\overline{\gamma}} \; , \label{e:coeff2}
\end{equation}
\begin{equation}
e^{2u} \hat{\Gamma}^\gamma_{\hat{0}\alpha} =
\Gamma^\gamma_{0\alpha} + 2 u_0 \delta_\alpha^\gamma +
i(\nabla^\gamma u_\alpha - 2 u_\alpha u^\gamma +
u^{\overline{\rho}} \Gamma^\gamma_{\overline{\rho}\alpha} - u^\rho
\Gamma^\gamma_{\rho\alpha} ) \; , \label{e:coeff3}
\end{equation}
where $u_A = T_A (u)$. Also $\nabla_{\overline{\beta}} u_\alpha =
T_{\overline{\beta}} (u_\alpha ) -
\Gamma^\mu_{\overline{\beta}\alpha} u_\mu$ and $\nabla^\gamma
u_\alpha = g^{\gamma\overline{\beta}} \nabla_{\overline{\beta}}
u_\alpha$, etc. For a proof of \eqref{e:coeff1}-\eqref{e:coeff3}
one may see Proposition 1 in \cite{kn:Dra2}, p. 39-40, or
\cite{kn:DrTo}, Chapter 2. One obtains
\[
{{\hat{R}\,_\alpha}^\beta}_{\lambda\overline{\mu}} =
{{R_\alpha}^\beta}_{\lambda\overline{\mu}} + 4 i
g_{\lambda\overline{\mu}} u_0 \delta_\alpha^\beta - 4(
\delta_\lambda^\beta g_{\overline{\mu}\alpha} +
\delta_\alpha^\beta g_{\overline{\mu}\lambda} ) u^\gamma u_\gamma
- \]
\[ - 2 g_{\alpha\overline{\mu}} \, \nabla_\lambda u^\beta
- 2 g_{\lambda\overline{\mu}} \nabla^\beta u_\alpha - 2
\delta^\beta_\alpha \nabla_{\overline{\mu}} u_\lambda - 2
\delta^\beta_\lambda \nabla_{\overline{\mu}} u_\alpha
\] hence (by the commutation formula $\nabla_{\overline{\beta}}
u_\alpha = \nabla_\alpha u_{\overline{\beta}} + 2 i
g_{\alpha\overline{\beta}} u_0$)
\[
e^{-2u} \hat{R}_{\alpha\overline{\beta}\lambda\overline{\mu}}
\xi^\alpha \xi^{\overline{\beta}} \xi^\lambda \xi^{\overline{\mu}}
= \]
\[ = R_{\alpha\overline{\beta}\lambda\overline{\mu}}
\xi^\alpha \xi^{\overline{\beta}} \xi^\lambda \xi^{\overline{\mu}}
- 8 \| \xi \|^4 \{ (\nabla_{\overline{\beta}} u_\alpha )
\eta^\alpha \eta^{\overline{\beta}} + u^\alpha u_\alpha - i u_0 \}
\] implying \eqref{e:relsect}.
\begin{corollary} Let $(M, \theta )$ be a compact Sasakian
manifold such that the Vaisman manifold $V = M \times S^1$ admits
an immersion $\phi$ into an ordinary complex Hopf manifold
$H_\Lambda$, $\Lambda = (\lambda , \cdots , \lambda )$, $0 <
\lambda  < 1$, and $\phi$ commutes with the Lee flows. Then
$H_\theta (\sigma ) \leq - \log \lambda$ for any $\sigma \in G_2
(H(M))$.
\end{corollary}
It should be observed that $\phi$ is obtained (cf. Theorem 5.1 in
\cite{kn:OrVe}) by first building an immersion $\tilde{V} \to H^0
(V^\prime , L^k_{\mathbb C})$ (cf. {\em op. cit.}, p. 139) of the
universal covering $\tilde{V}$ into a suitable space of
holomorphic sections and the problem of the effective
computability of $n$ and $\Lambda$ (in terms of the given data,
i.e. the locally conformal K\"ahler structure on $V$) is an open
problem.
\begin{corollary} Let $u$ be a CR-pluriharmonic function on $M$
i.e. there is a $C^\infty$ function $v : M \to {\mathbb R}$ such
that $u + iv$ is a CR function. Then $e^{2u} H_{\hat{\theta}}
(\sigma ) \leq H_\theta (\sigma ) + 2 v_0$ for any $\sigma \in G_2
(H(M))$.
\end{corollary}
{\em Proof}. By a result of J. M. Lee, \cite{kn:Lee}, if $u$ is
CR-pluriharmonic and $v$ is conjugate to $u$ then the complex
Hessian of $u$ is given by
\[ \nabla_{\overline{\beta}} u_\alpha = (i u_0 - v_0 )
g_{\alpha\overline{\beta}} \] hence (by \eqref{e:relsect}) $e^{2u}
H_{\hat{\theta}} (\sigma ) = H_\theta (\sigma ) - 2 u^\alpha
u_\alpha + 2 v_0 \leq H_\theta (\sigma ) + 2 v_0$. Q.e.d.
\begin{appendix}
\section{The Sasakian Schur theorem}
As a first purpose of Appendix A we prove Theorem \ref{t:4}. First
note that (by Proposition 7.2 in \cite{kn:KoNo}, Vol. II, p. 167)
the $4$-tensor \eqref{e:21} satisfies
\begin{equation}
\label{e:A}
R_0 (X,Y,Z,W) = - R_0 (Y,X, Z,W) = - R_0 (X,Y,W,Z) \; ,
\end{equation}
\begin{equation}
\label{e:B} R_0 (X,Y,Z,W) = R_0 (Z,W, X,Y) \; ,
\end{equation}
\begin{equation}\label{e:C} \sum_{YZW} R_0 (X,Y,Z,W) = 0 \; ,
\end{equation}
\begin{equation}\label{e:D}
R_0 (J X , J Y , Z , W) = R_0 (X,Y, JZ , J W) = R_0 (X,Y,Z,W) \; ,
\end{equation}
for any $X,Y,Z,W \in H(M)$. On the other hand (cf. e.g. the
Appendix A in \cite{kn:BaDr})
\begin{equation}
\label{e:A1} R(X,Y,Z,W) = - R(Y,X,Z,W) = - R(X,Y,W,Z) \; ,
\end{equation}

\begin{equation}
\label{e:B1} R(X,Y,Z,W) = R(Z,W, X,Y) -
\end{equation}
\[  -
2 \Omega (Y,Z) A(X,W) + 2 \Omega (Y,W) A(X,Z) - \]
\[ - 2\Omega (X,W) A(Y,Z) + 2 \Omega (X,Z) A(Y,W) \; ,
\]
\begin{equation}
\label{e:C1}  \sum_{YZW} R(X,Y,Z,W) = - 2 \sum_{YZW} \Omega (Y,Z)
A(W,X) \; ,
\end{equation}
for any $X,Y,Z,W \in H(M)$. Moreover (by $\nabla J = 0$)
\begin{equation}
\label{e:D1.i} R(JX, J Y , Z, W) = R(X,Y,Z,W) \; .
\end{equation}
Let us assume that $H_\theta = f \circ \pi$ for some $f : M \to
{\mathbb R}$ and set $L = R - 4 f R_0$. The properties
\eqref{e:A}-\eqref{e:D} and \eqref{e:A1}-\eqref{e:D1.i} imply
\begin{equation}
\label{e:A2}
L(X,Y,Z,W) = - L(Y,X, Z,W) = - L(X,Y,W,Z) \; ,
\end{equation}
\begin{equation}
\label{e:B2} L(X,Y,Z,W) = L(Z,W,X,Y) - \end{equation}
\[  -
2 \Omega (Y,Z) A(X,W) + 2 \Omega (Y,W) A(X,Z) - \] \[ - 2 \Omega
(X,W) A(Y,Z) + 2 \Omega (X,Z) A(Y,W) \; ,
\]
\begin{equation} \label{e:C2} \sum_{YZW} L(X,Y,Z,W) = - 2
\sum_{YZW} \Omega (Y,Z) A(W,X) \; ,
\end{equation}
\begin{equation} \label{e:D2.i} L(JX, J Y, Z,W) = L(X,Y,Z,W) \; .
\end{equation}
As to the analog of the second equality in \eqref{e:D} for the
$4$-tensor $L$ one has (by \eqref{e:B2} and \eqref{e:D2.i})
\[
L(X,Y, J Z , J W) = L(J Z , J W , X , Y) - \]
\[
- 2 \Omega (Y, JZ) A(X, J W) + 2 \Omega (Y, J W) A(X , J Z) -
\]
\[ - 2 \Omega (X , J W) A(Y, J Z) + 2 \Omega (X, J Z) A(Y, J
W) =
\]
\[
= L(Z,W,X,Y) + \]
\[ + 2 g_\theta (Y,Z) A(X , J W) - 2 g_\theta
(Y,W) A(X , J Z) +
\]
\[
+ 2 g_\theta (X,W) A(Y , J Z) - 2 g_\theta (X,Z) A(Y, J W)
\] and applying once more \eqref{e:B2} leads to
\begin{equation}
\label{e:D2.ii}
L(X,Y, J Z, J W) = L(X,Y,Z,W) + \end{equation}
\[
+ 2 \Omega (Y,Z) A(X,W) - 2 \Omega (Y,W) A(X,Z)  +\]
\[+ 2 \Omega
(X,W) A(Y,Z) - 2 \Omega (X,Z) A(Y,W) + \]
\[ + 2 g_\theta (Y,Z) A(X, J W) - 2
g_\theta (Y,W) A(X, J Z) +\] \[   + 2 g_\theta (X,W) A(Y, J Z) - 2
g_\theta (X,Z) A(Y, J W) \; .
\]
\end{appendix}
Note that (by the very definition of $H_\theta$)
\begin{equation}
\label{e:1.app} L(X, J X , X , J X) = 0
\end{equation}
for any $X \in H(M)$. We shall need the $4$-tensor $K$ defined by
\[
K(X,Y,Z,W) = L(X, J Y , Z , J W) + \]
\[  + L(X , J Z , Y , J
W) + L(X , J W , Y , J Z) \; . \] As an immediate consequence of
\eqref{e:1.app}
\begin{equation}
K(X,X,X,X) = 0 \; . \label{e:36}
\end{equation}
Using repeatedly \eqref{e:A2}-\eqref{e:B2} and
\eqref{e:D2.i}-\eqref{e:D2.ii} together with
\[ \Omega (J X , J Y) = \Omega (X,Y), \;\;\; A(J X , J Y ) = -
A(X,Y) \; , \] we derive the identities
\begin{equation}
\label{e:2.app} K(Y,X,Z,W) = K(X,Y,Z,W) + 4 \Omega (X,Y) A (Z,W) -
\end{equation}
\[ - 2 g_\theta (X,Z) A(Y, J W) + 2 g_\theta
(Y,W) A(X , J Z) - \]
\[  - 2 g_\theta (X,W) A(Y, J Z) + 2
g_\theta (Y,Z) A(X , J W) \; ,
\]

\begin{equation}
\label{e:3.app}
K(X,Y,W,Z) = K(X,Y,Z,W) +
\end{equation}
\[ +
2 g_\theta (Y,Z) A(X, J W) - 2 \Omega (Y,W) A(X,Z) - \]
\[ - 2
g_\theta (X,W) A(Y, J Z) + 2 \Omega (X,Z) A(Y,W) + \]
\[ + 2
\Omega (Y,Z) A(X,W) - 2 g_\theta (Y,W) A(X, J Z) - \]
\[  - 2
\Omega (X,W) A(Y,Z) + 2 g_\theta (X,Z) A(Y, J W) \; ,
\]

\begin{equation}
\label{e:4.app} K(Z,Y,X,W) = K(X,Y,Z,W) + \end{equation}
\[ + 2 \Omega (Z,W) A(X,Y) + 2 \Omega (X,Y) A(Z,W) + \]
\[ + 2 \Omega (W,X) A(Y,Z) + 2 \Omega (Y,Z) A(X,W) \; ,
\]

\begin{equation}
\label{e:5.app} K(X, W, Z, Y) = K(X, Y, Z, W) + 4 \Omega (X,Y)
A(Z,W) - \end{equation}
\[ - 4\Omega (Y,W) A(X,Z) - 4 \Omega
(X,W) A(Y,Z) + \]
\[  + 4 g_\theta (X,Y) A(Z, J W)
- 4 g_\theta (Z,W) A(X, J Y) + \]
\[ + 4 g_\theta (Y,Z) A (X, J W) - 4 g_\theta (X,W) A(Y, J Z) \; ,
\]

\begin{equation}
\label{e:6.app} K(W, Y, Z, X) = K(X,Y, Z, W) - 4
\Omega (X, W) A(Y,Z) +
\end{equation}
\[ + 2 g_\theta (X,Y) A(Z, J W) + 2 g_\theta (X,Z) A(Y,
J W) - \]
\[ - 2 g_\theta (Y,W) A(X, J Z) - 2 g_\theta (Z, W)
A(X, J Y) \; ,\]

\begin{equation}
\label{e:7.app}
K(X,Z, Y,W) = K(X,Y,Z,W) +
\end{equation}
\[ +
2 g_\theta (Y,W) A(X, J Z) + 2 \Omega (Z,W) A(X,Y) - \]
\[ - 2 g_\theta (X,Z) A(Y, J W) + 2 \Omega (X,Y) A(Z,W)
- \]
\[- 2 \Omega (Y,W) A(X,Z) - 2 g_\theta (Z,W) A(X, J Y)
- \]
\[- 2 \Omega (X,Z) A(Y,W) + 2 g_\theta (X,Y) A(Z, J W) \; . \]
Let us replace $X$ by $X+Y$ in \eqref{e:36} and use
\eqref{e:2.app}-\eqref{e:7.app})
 so that to obtain
\begin{equation}
\label{e:43} 2 K(X,X,X,Y) + 3 K(X,Y,X,Y) + 2 K(X, Y, Y, Y) =
\end{equation}
\[  = 8 \{ g_\theta (X,Y) A(X, J X) - g_\theta
(X,X) A(X, J Y)  + \Omega (X,Y) A(X,X) \} \; .
\]
Next we replace $Y$ by $Y + Z$ in
\eqref{e:43} and use again \eqref{e:2.app}-\eqref{e:7.app}. We
obtain
\begin{equation}
\label{e:44}
K(X,Y,X,Z) + K(X,Y,Y,Z) + K(X,Y,Z,Z) =
\end{equation}
\[
= 2 \Omega (Y,Z) A(X,X) - 2 \Omega (X,Y) A(X,Z) + 2 \Omega (X,Z)
A(X,Y) - \] \[ - 2 \Omega (X,Y) A(Y,Z)  + 2 \Omega (Y,Z) A(X,Y) -
2 \Omega (X,Y) A(Z,Z) +\]
\[
+ 2 \Omega (Y,Z) A(X,Z) + 2 \Omega (X,Z) A(Y,Y) + 2 \Omega (X,Z)
A(Y,Z) + \]
\[ + 4 g_\theta (X,Z) A(X, J Y) - 4 g_\theta (X,Y)
A(X, J Z) +  \]
\[  + 2 g_\theta (X,Z) A(Y, J Z) - 2
g_\theta (X,Y) A(Y, J Z) + \]
\[ + 2 g_\theta
(X,Z) A(Y, J Y) - 2 g_\theta (X,Y) A(Z, J Z) + \]
\[ + 2 g_\theta (Y,Z) A(X, J Y) - 2 g_\theta (Y,Z) A(X,
J Z) + \]
\[  + 2 g_\theta (Z,Z) A(X, J Y) - 2 g_\theta (Y,Y)
A(X, J Z) \; . \] Finally let us replace $Z$ by $Z + W$ in
\eqref{e:44} and derive the expression of the $4$-tensor $K$
\begin{equation}
\label{e:45} L(X,J Y, Z, J W) + L(X, J Z, Y, J W) + L(X, J W, Y,
J Z) =
\end{equation}
\[ = 2 \Omega (Y,W) A(X,Z) + 2 \Omega
(X,W) A(Y,Z) + 2 \Omega (Y,X) A(Z,W)  +\] \[ + 2 g_\theta (X,W)
A(Y, J Z) - 2 g_\theta (Y,Z) A(X, J W) + \]
\[ + 2 g_\theta (Z,W) A(X, J Y) - 2 g_\theta (X,Y) A(Z,
J W) \; . \] Setting $Z = X$ and $W = Y$ in \eqref{e:45} gives
\[ L(X, JY, X, J Y) + L(X, J X , Y, J Y) - L(X, J Y, J X , Y) = 0 \; .
\] Let us apply \eqref{e:D2.ii} to the last term. We obtain
\begin{equation}
\label{e:*} 2 L(X, J Y, X, J Y) + L(X, J X, Y, J Y) = 4 \Omega
(X,Y) A(X,Y) +
\end{equation}
\[  + 2 g_\theta (Y,Y) A(X, J X)
- 2 g_\theta (X,X) A(Y, J Y) \; . \] On the other hand we replace
$(Y,Z,W)$ in \eqref{e:C2} by $(J X, Y, J Z)$ so that to get
\[
L(X, J X, Y, J Y) + L(X,Y, J Y, J X) + L(X, J Y, J X, Y) =
\]
\[ = - 2 \{ \Omega (J X, Y) A(J Y , X) + \Omega (Y, J Y) A(J
X , X) + \Omega (J Y, J X) A(Y,X) \} \] or (again by
\eqref{e:D2.ii})
\begin{equation}
\label{e:**} L(X, JX, Y, J Y) - L(X,Y,X,Y) - L(X, J Y, X , J Y)
= \end{equation}
\[ = 2 g_\theta (X,Y) A(X, J Y) - 2 g_\theta
(Y,Y) A(X, J X)  - 2 \Omega (X,Y) A(X,Y) \; .\] Let us subtract
\eqref{e:**} from \eqref{e:*}. This gives
\begin{equation}
\label{e:alpha}
3 L(X, J Y, X, J Y) + L(X,Y,X,Y) =
\end{equation}
\[ = 6 \Omega (X,Y) A(X,Y) - 2 g_\theta (X,X) A(Y, J
Y) + \]
\[ + 4 g_\theta (Y,Y) A(X, J X) - 2 g_\theta (X,Y)
A(X, J Y) \; .\] Replacing $Y$ by $J Y$ in \eqref{e:alpha} leads
to the identity
\begin{equation}
\label{e:beta}
3 L(X,Y,X,Y) + L(X, JY, X, JY) =
\end{equation}
\[ = - 6 g_\theta (X,Y) A(X, J Y) + 2 g_\theta (X,X)
A(Y, J Y) + \]
\[  + 4 g_\theta (Y,Y) A(X, J X) + 2 \Omega
(X,Y) A(X,Y) \; . \] Solving for $L(X,Y,X,Y)$ in the linear system
\eqref{e:alpha}-\eqref{e:beta} gives
\begin{equation}
\label{e:50}
L(X,Y,X,Y) = g_\theta (X,X) A(Y, J Y) -
\end{equation}
\[ - 2 g_\theta (X,Y) A(X, J Y) + g_\theta (Y,Y) A(X, J X) \; .
\]
Once $L(X,Y,X,Y)$ is known one may apply \eqref{e:A2}-\eqref{e:C2}
and the algebraic scheme in the proof of Proposition 1.2 in
\cite{kn:KoNo}, Vol. I, p. 198, to compute the whole of
$L(X,Y,Z,W)$. Precisely we replace $Y$ by $Y+W$ in \eqref{e:50}
and use \eqref{e:C2} so that to get
\[ L(X,Y,X,W) = \Omega (X,Y) A(X,W) + \]
\[ + \Omega (Y,W) A(X,X) - \Omega (X,W) A(X,Y) + \]
\[+ g_\theta (X,X) A(Y, J W) - g_\theta (X,Y) A(X, J W) - \]
\[ - g_\theta (X, W) A(X, J Y) + g_\theta (Y,W) A(X, J X) \; .
\] Next we replace $X$ by $X + Z$ and derive the identity
\begin{equation}
\label{e:51}
L(X,Y,Z,W) = L(X,W,Y,Z) + 2 \Omega (X,Z) A(Y,W) -
\end{equation}
\[ - \Omega (X,Y) A(Z,W) + \Omega (Z,W) A(X,Y)
- \]
\[ - \Omega (Y,Z) A(X,W) + \Omega (W,X)
A(Y,Z) + \]
\[ + 2 g_\theta (X,Z) A(Y, J W) + 2
g_\theta (Y,W) A(X, J Z) - \]
\[ - g_\theta (X,Y)
A(Z, J W) - g_\theta (Y,Z) A(X, JW) - \]
\[ - g_\theta (X,W)
A(Y, J Z) - g_\theta (Z,W) A(X, J Y) \; .
\]
Another identity of the kind is got by replacing $(Y,Z,W)$ in
\eqref{e:51} by $(Z,W,Y)$ i.e.
\begin{equation}
\label{e:52}
L(X,Z,W,Y) =
L(X,Y,Z,W) + 2 \Omega (X,W) A(Z,Y) -
\end{equation}
\[ -
\Omega (X,Z) A(W,Y) + \Omega (W,Y) A(X,Z) - \]
\[ - \Omega (Z,W) A(X,Y) + \Omega (Y,X) A(Z,W) + \]
\[ + 2 g_\theta (X,W) A(Z, J Y) + 2 g_\theta (Z,Y)
A(X, J W) - \]
\[  - g_\theta (X,Z) A(W, J Y) -
g_\theta (Z,W) A(X, J Y) - \]
\[  - g_\theta (X,Y) A(Z, J W) -
g_\theta (W, Y) A(X, J Z) \; .\] Finally let us compute $3
L(X,Y,Z,W)$ by expressing the second and third copy of
$L(X,Y,Z,W)$ from \eqref{e:51}-\eqref{e:52}, respectively. Then
(by \eqref{e:C2})
\[ L(X,Y,Z,W) = \Omega (X,Z) A(Y,W) - \Omega (Y,Z) A(X,W) + \]
\[ + \Omega (Y,W) A(X,Z) - \Omega (X,W) A(Y,Z) + \]
\[ + g_\theta (X,Z) A(Y, J W) - g_\theta (Y,Z) A(X, J W) + \]
\[ + g_\theta (Y,W) A(X, J Z) - g_\theta (X,W) A(Y, J Z)
\] and \eqref{e:19} in Theorem \ref{t:4} is proved. The identity
\eqref{e:20} follows from \eqref{e:19} by contraction. Another
scope of Appendix A is to establish the following Sasakian analog
to the complex Schur theorem in \cite{kn:KoNo}, Vol. II, p. 168.
\begin{theorem} Let $(M, \theta )$ be a connected pseudohermitian manifold
of CR dimension $n \geq 3$. Assume that $H_\theta = f \circ \pi$
for some $C^\infty$ function $f : M \to {\mathbb R}$. If $S = 0$
then $\nabla f = \theta (\nabla f) T$. Moreover if $\tau = 0$ then
$f$ is constant. \label{t:5}
\end{theorem}
Therefore each Sasakian manifold of CR dimension $\geq 3$ whose
pseudohermitian sectional curvature \eqref{e:1} is but a point
function is actually a pseudohermitian space form. The proof of
the complex Schur theorem is to show that each K\"ahlerian
manifold whose holomorphic sectional curvature is a point function
$f$ is an Einstein manifold. Yet, if the manifold dimension is
$\geq 3$, the Einstein condition together with the second Bianchi
identity imply that $f$ is constant (cf. Note 3 in \cite{kn:KoNo},
Vol. I, p. 292-294). As argued in \cite{kn:Lee}, the
pseudo-Einstein condition together with the second Bianchi
identity (for the Tanaka-Webster connection) does not imply in
general that the pseudohermitian scalar curvature is constant (due
to the presence of torsion terms in the second Bianchi identity).
Therefore we use the full curvature tensor
\[ R(X,Y,Z,W) = f \{ 2 \Omega (X,Y) \Omega (Z,W) + \]
\[ + g_\theta (X,Z) g_\theta (Y,W) - g_\theta (X,W) g_\theta (Y,Z) +\]
\[ + \Omega (X,Z) \Omega (Y,W) - \Omega (X,W) \Omega (Y,Z)\} + \]
\[+ g_\theta (X,Z) A(Y, J W) - g_\theta (X,W) A(Y, J Z) + \]
\[ + g_\theta (Y,W) A(X, J Z) - g_\theta (Y,Z) A(X, J W) + \]
\[ + \Omega (X,Z) A(Y,W) - \Omega (X,W) A(Y,Z) + \]
\[ + \Omega (Y,W) A(X,Z) - \Omega (Y,Z) A(X,W) \] for any $X,Y,Z,W
\in H(M)$. Rather than contracting we take the covariant
derivative of the previous identity. A rather lengthy calculation
(based on $\nabla g_\theta = 0$ and $\nabla \Omega = 0$) leads to
\begin{equation}
\label{e:53}
(\nabla_U R)(X,Y,Z,W) = U(f) \{ 2 \Omega (X,Y)
\Omega (Z,W) +
\end{equation}
\[ + g_\theta (X,Z) g_\theta
(Y,Z) - g_\theta (X,W) g_\theta (Y,Z) +  \]
\[  +\Omega (X,Z) \Omega (Y,W) - \Omega (X,W) \Omega (Y,Z) \} + \]
\[ + g_\theta (X,Z) (\nabla_U A)(Y, J W) - g_\theta
(X,W) (\nabla_U A)(Y, J Z) + \]
\[  + g_\theta (Y,W) (\nabla_U A)(X, J Z) - g_\theta (Y,Z) (\nabla_U A)(X, J W)
+ \]
\[ + \Omega (X,Z) (\nabla_U A)(Y,W) - \Omega
(X,W) (\nabla_U A)(Y,Z) + \]
\[  + \Omega (Y,W) (\nabla_U
A)(X,Z) - \Omega (Y,Z) (\nabla_U A)(X,W) \] for any $U \in T(M)$
and any $X,Y,Z,W \in H(M)$. From now on let $U \in H(M)$ and let
us take the cyclic permutation $U \to Z \to W \to U$ in
\eqref{e:53} to get two more identities of the kind. Next let us
add up the three resulting identities and express $\sum_{UZW}
(\nabla_U R)(X,Y,Z,W)$ from the second Bianchi identity
\begin{equation}
(\nabla_U R)(X,Y,Z,W) + (\nabla_Z R)(X,Y,W,U)
+ (\nabla_W R)(X,Y,U,Z) =
\end{equation}
\[ = g_\theta (\Omega
(U,Z) W + \Omega (Z,W) U + \Omega (W,U) Z \, , \,  S(Y,X))
\] for any $X,Y,Z,W,U \in H(M)$. A calculation based on
\[ (\nabla_X A)(Y,Z) - (\nabla_Y A)(X, Z) = g_\theta (S(X,Y), Z) \; ,
\]
\[ (\nabla_X A)(J Y , Z) - (\nabla_Y A)(J X , Z) = g_\theta
(S(X,Y), J Z) \; , \] furnishes
\begin{equation}
\label{e:55} g_\theta ( \Omega (U,Z) W + \Omega (Z,W) U + \Omega
(W,U) Z \, , \, S(Y, \cdot ))^\sharp =
\end{equation}
\[ = U(f) \{ g_\theta (Y,W) Z - g_\theta (Y,Z) W
 + 2 \Omega (Z,W) J Y + \]
 \[ +\Omega (Y,W) J Z - \Omega
(Y,Z) J W \} + \]
\[ + Z(f) \{ g_\theta (Y,U) W -
g_\theta (Y,W) U  + 2 \Omega (W,U) J Y +\] \[ + \Omega (Y,U) J W -
\Omega (Y,W) J U \} + \]
\[ +
W(f) \{ g_\theta (Y,Z) U - g_\theta (Y,U) Z + 2 \Omega (U,Z) J Y
+\] \[ + \Omega (Y,Z) J U - \Omega (Y,U) J Z \} + \]
\[ + g_\theta (S(U,W), J Y) Z -
g_\theta (S(U,Z), J Y) W  - g_\theta (S(Z,W) , J Y) U + \]
\[ + g_\theta (Y,Z) J \, S(U,W) - g_\theta (Y,W) J \,
S(U,Z)  - g_\theta (Y,U) J \, S(Z,W) + \]
\[ + g_\theta (S(U,W), Y) J Z - g_\theta (S(U,Z), Y) J W - g_\theta
(S(Z,W), Y) J U + \]
\[ + \Omega (Y,W) S(U,Z) - \Omega (Y,Z)
S(U,W) + \Omega (Y,U) S(Z,W) \; .\] Let $U \in H(M)$ be arbitrary
and let us choose $Y, Z \in H(M)$ such that $g_\theta (Y,Z) =
\Omega (Y, Z) = 0$, $g_\theta (Y,U) = \Omega (Y,U) = 0$ and
$g_\theta (Z, U ) = \Omega (Z,U) = 0$. Also we assume $\| Y \| =
1$ and set $W = Y$. Then \eqref{e:55} gives
\[
U(f) Z - Z(f) U - J \, S(U,Z) + \]
\[ + \Omega (S(Z,U), Y) Y + g_\theta (S(Z,U), Y) J Y
+ \Omega (S(U,Y), Y) Z +\]
\[ + g_\theta (S(U,Y), Y) J Z
+ \Omega (S(Y,Z), Y) U + g_\theta (S(Y,Z), Y) J U = 0 \; . \]
Hence $S = 0$ yields $U(f) = 0$ and the first statement in Theorem
\ref{t:5} is proved. Note that (as $S = 0$)
\[ (\nabla_Z R)(X,Y,W,T) = 0 \; . \]
Hence, by the second Bianchi identity
\[ (\nabla_T R)(X,Y,Z,W) + (\nabla_Z R)(X,Y,W,T) + (\nabla_W
R)(X,Y,T,Z) = \]
\[= g_\theta (R(\tau (W), Z) Y , X) -
g_\theta (R(\tau (Z), W) Y , X) -\] \[ - g_\theta (R(T_\nabla
(Z,W), T)Y , X) \] we obtain
\begin{equation}
\label{e:56} (\nabla_T R)(X,Y,Z,W) =
\end{equation}
\[ = g_\theta (R(\tau (W), Z) Y \, , \, X) - g_\theta (R(\tau (Z), W)
Y \, , \, X)\] for any $X,Y,Z,W \in H(M)$. Let us set $U = T$ in
\eqref{e:53} and substitute from \eqref{e:56} in the resulting
identity. We obtain
\begin{equation}
\label{e:57} R(\tau (W), Z) Y - R(\tau (Z), W) Y = T(f) \{
g_\theta (Y,W) Z - g_\theta (Y,Z) W +
\end{equation}
\[  + 2 \Omega (Z,W) J Y + \Omega (Y,W) J Z -
\Omega (Y,Z) J W \} + \]
\[  + (\nabla_T A)(Y, J
W) Z - (\nabla_T A)(Y, J Z) W  + g_\theta (Y,W) (\nabla_T \tau ) J
Z -\] \[ - g_\theta (Y,Z) (\nabla_T \tau ) J W  + (\nabla_T
A)(Y,W) J Z - (\nabla_T A)(Y,Z) J W   +\] \[ + \Omega (Y,W)
(\nabla_T \tau )Z - \Omega (Y,Z) (\nabla_T \tau )W \; .
\] Let $Y,Z \in H(M)$ such that $\| Y \| = 1$ and $g_\theta (Y,Z) =
\Omega (Y,Z) = 0$ and set $W = Y$ in \eqref{e:57}. Together with
the assumption $\tau = 0$ this yields $T(f) = 0$. Theorem
\ref{t:5} is proved.

\end{document}